\documentclass{article}
\usepackage{graphicx}
\usepackage{amssymb}
\usepackage{amsmath}
\usepackage{a4wide}
\usepackage{caption}
\usepackage{subcaption}
\usepackage[affil-it]{authblk}
\usepackage{xr}
\usepackage{color}
\usepackage[utf8]{inputenc}
\usepackage[T1]{fontenc}
\newtheorem{theorem}{Theorem}[section]

\newtheorem{corollary}[theorem]{Corollary}

\newtheorem{definition}[theorem]{Definition}

\newtheorem{remark}[theorem]{Remark}

\begin{document}

\title{Application of planar Randers geodesics \\ with river-type perturbation in search models}
\author{Piotr Kopacz}
\date{\texttt{}}
\maketitle
\begin{abstract}
\noindent
We consider remodeling the planar search patterns, in the presence of the river-type perturbations represented by the weak vector field, basing on the time-optimal paths as Finslerian solutions to the Zermelo navigation problem via Randers metric.
\end{abstract}
\smallskip
\textbf{M.S.C. 2010}: 53C22, 53C60, 53C21, 53B20, 68U35, 68T20.
\smallskip


\smallskip


\section{Introduction}

\subsection{Motivation}

In geometric geodesy there are two standard problems. First (direct) geodetic problem is to determine the final point $A_2 \in M$ represented by its coordinates, when the starting point $A_1\in M$ is given as well as the distance and the initial angle of the geodesic of modeling surface $M$ coming from $A_1$ towards $A_2$, in reference to a fixed direction. Second (inverse) geodetic problem is to determine the length of the geodesic between the endpoints and its angles at them, when both points are given. More advanced, what covers the area of our interest, these problems are applied in navigation in the complex route planning and monitoring processes under the real perturbations resulting in leeway and drift. Acting perturbation can modify the preplanned trajectory and change the endpoint in the former problem in each part of the passage. Alternatively, the correction is applied to keep the ship on set track over ground. Thus, to determine the correction as a function of position and time in general, we need to find the notions considered in the latter problem. In the paper we aim to investigate the special solutions to generalized geodetic problems, which additionally provide the time-optimality under acting perturbation modeled by a vector field. We achieve the optimality making use of Randers metric. Hence, the modern approach to Zermelo's problem (1931) \cite{zermelo} in Finsler geometry is considered (cf. \cite{colleen_shen}). Due to potential real applications we pay more attention to low dimensions in our research (cf. \cite{kopi}).  


\subsection{Preliminaries and posing the problem}

Let a pair $ (M,h) $ be a Riemannian manifold where $h = h_{ij}dx^i\otimes dx^j$ is a Riemannian metric  and denote the corresponding norm-squared of tangent vectors $\mathbf{y} \in T_x M$ by
$\left|\mathbf{y} \right|^2 = h_{ij}y^iy^j = h(\mathbf{y}, \mathbf{y}).$
The unit tangent sphere in each $T_x M$ consists of all tangent vectors $u$ such that $\left|u\right| = 1$. Equivalently, $h(u, u) = 1$. Then we introduce a vector field $W$ such that $\left|W \right| < 1$, thought of as the spatial velocity vector of a weak wind on the Riemannian sea $(M, h)$. Before $W$ sets in, a passage from the base to the tip of any $u$ would take one unit of time. The effect of the wind is to cause the journey to veer off course or just off target if $u$ is collinear with $W$. Instead of $u$ we traverse the resultant $v = u + W$ within the same one unit of time. Thus, in the presence of the wind the Riemannian metric $h$ no longer gives the travel time along vectors. This causes the introduction of a function $F$ on the tangent bundle $T M$, in order to keep track of the travel time needed to traverse tangent vectors $\mathbf{y}$ under perturbation. For all those resultants $v = u + W$, we have $F(v) = 1$. Within each tangent space $T_x M$, the unit sphere of $F$ is the $W$-translate of the unit sphere of $h$. Since this $W$-translate is no longer centrally symmetric, $F$ cannot possibly be Riemannian \cite{colleen_shen}. 

In general, the problem is to find the trajectory followed by the ship and the corresponding steering angle such that the ship completes her journey in the least time. In the paper we focus on 2D river-type perturbation, i.e with one component zeroed. Setting in the planar coordinate system let the vertical to be the one. Now we consider a ship navigating on a sea with a current and assume that the current runs parallel to the $x$-axis. The state of a ship is represented by its position $(x,y)$ and heading angle $\varphi$. Thus, the maximal state space is $\mathbb{R}^2 \times S^1$ if the convexity condition is fulfilled in the whole plane. Otherwise we have to restrict our domain for the position coordinates as a subset of $\mathbb{R}^2$. The current may have a rotational effect as well as a translational effect on the ship, and this effect depends on the ship's heading angle. Note that in the general Zermelo navigation problem we consider the time-dependent vector field which is given in the most general form 
\begin{equation}
\label{general_field}
W = \frac{\partial}{\partial t} + W^i\left(t,x^j\right)\frac{\partial}{\partial x^i}.
\end{equation}
\noindent
The variational approach to the solution to the Zermelo navigation problem in time-space under time-dependent perturbation (\ref{general_field}) via the Euler-Lagrange equations is presented in \cite{palacek} and followed in \cite{kopi}. In the paper let us begin with the stationary current, namely the perturbation does not depend on time. This is a special case of the Zermelo navigation problem researched by C. Caratheodory by means of Hamiltonian formalism in the calculus of variations (cf. \S 276 in \cite{caratheodory}). We aim to find the deviation of calm sea geodesics under the action of distributed wind which is modeled by the vector field on manifold $M$. We aim to find the deviation of the Riemannian geodesics with application of special class of Finsler metrics, namely Randers metrics. Consequently, this also enables us to analyze the geometric properties of the trajectories as Finslerian solutions to Zermelo's problem depending on the type of perturbation and the initial or boundary conditions. 


\section{Randers metric}
From \cite{colleen_shen} we know that in Riemann-Finsler geometry Randers metrics may be identified with the solutions to the navigation problem on Riemannian manifolds. This navigation structure establishes a bijection between Randers spaces and pairs $(h,W)$ of Riemannian metrics $h$ and vector fields $W$ on the manifold $M$. 


\subsection{The general form of the metric}

The resulting Randers metric is composed of the new Riemannian metric and $1$-form, and is given by
\begin{equation}
\label{ran}
F(\mathbf{y}) = \frac{ \sqrt{ \left[ h(W,\mathbf{y}) \right]^2 + |\mathbf{y}|^2 \lambda} } {\lambda} - \frac{ h(W,\mathbf{y})} {\lambda}.
\end{equation}

\noindent
The resulting Randers metric can also be presented in the form $F = \alpha + \beta$ as the sum of two components. Explicitly,
\begin{itemize}
\item  the first term is the norm of $\mathbf{y}$ with respect to a new Riemannian metric
\begin{equation}
\alpha(x,\mathbf{y})  = \sqrt{a_{ij}(x)y^iy^j}, \quad \text{ where } a_{ij} = \frac{ h_{ij}}{\lambda} + \frac{W_i}{\lambda}\frac{W_j}{\lambda},
\label{Alpha}
\end{equation}

\item the second term is the value on $\mathbf{y}$ of a differential 1-form
\begin{equation}
\beta(x,\mathbf{y})=b_i(x)y^i,\quad \text{where}\quad  b_i = \frac{-W_i}{\lambda}.
\label{Beta}
\end{equation}
\end{itemize}
where $W_i=h_{ij}W^j$ and $\lambda=1-W^iW_i$. 

As we aim to combine the theoretical research with real applications we have to make a remark here in reference to the notion of perturbation. In marine navigation there are two main types of real perturbation considered, namely wind and current (stream). Both of them affect the ship's trajectories and the corresponding angles, i.e. heading, course through the water and course over ground. In the literature on Finsler geometry the notions of wind and current are treated subsequently and equivalently by different authors. We will also use both notions of acting perturbation subsequently, although the notion of current fits better to the river-type perturbation which we consider in the paper. However, we let ourselves to call acting vector field $W$ equivalently as a wind in order to be in accordance with the contributions to the navigation problem in Finslerian approach. To be more precise, note that the notion of current or equivalently wind, which are used by different authors, means in fact the total drift of a vessel, i.e. the resulting perturbation, and so it ought to be considered in other papers implicitly if not otherwise stipulated.


\subsection{Randers metric with Euclidean background and river-type perturbation}

Let us assume that the initial Riemannian metric $h_{ij}$ to be perturbed is the standard Euclidean metric $\delta_{ij}$ on $\mathbb{R}^2$ as it can be applied in local modeling referring to Zermelo's problem. The Randers metric $\eqref{ran}$ which comes from the navigaton data $(h, W)$ including the initial Euclidean metric is expressed as follows 
\begin{equation}
F(x,\mathbf{y}) = \frac{ \sqrt{( \delta_{ij}W^iy^j)^2 + |\mathbf{y}|^2 (1-|W|^2)} } {1-|W|^2} - \frac{ \delta_{ij}W^iy^j} {1-|W|^2}.
\end{equation}
Since our focus is on dimension two, we denote the position coordinates $(x^1, x^2)$  by $(x, y)$, and expand arbitrary tangent vectors $y^1\frac{\partial}{\partial x^1}+y^2\frac{\partial}{\partial x^2}$ at  $(x^1, x^2)$ as $(x,y;u,v)$ or $u\frac{\partial}{\partial x} + v\frac{\partial}{\partial y}$. Thus, adopting the notations in two dimensional case yields
\begin{equation}
  \label{W1W2}
F(x,y; u,v) = \frac{\sqrt{u^2+v^2-(uW^2-vW^1)^2}-W^1u-W^2v}{1-|W|^2},
\end{equation}
\noindent
what defines the metric for an arbitrary vector field $W=(W^1, W^2)$. Let us consider the perturbation $W_I$ represented by the river-type vector field $W$ given in the following form
\begin{alignat}{1}
\label{pole_river}
(W)&	\left\{ \begin{array}{l l} 
		W^1 = W^1(x,y) \\
		W^2 = W^2(x,y) 	
	\end{array} \right.
	\xrightarrow{}
	(W_I)\left\{ \begin{array}{l l} 
		W^1(x,y):= f(y) \\
		W^2(x,y) := 0
	\end{array} \right.
\end{alignat}
what determines the scenario for the navigation problem. To apply Finslerian version of the Zermelo navigation we have to take into consideration the assumptions which limit the set of the solutions. Namely, only weak perturbation can be applied, i.e. $|W|<1$ and vessel's own speed is constant, $|(u, v)|=1$. The former assumption is due to the requirement of the strong convexity of the Finsler metric. This is implied by Theorem \ref{THM} which we shall apply next. Let us also make two remarks noting Finslerian research on the problem under stronger perturbation.    
\begin{remark}
For stronger perturbation, i.e. $|W|=1$, Kropina metric has been applied to the Zermelo navigation problem in Finsler geometry recently (R.Yoshikawa, S. Sabau in \cite{kropina}). 
\end{remark}
\begin{remark}
Strong wind, i.e. $|W|>1$, in the Zermelo navigation has been considered very recently (e.g. E. Caponio, M.A. Javaloyes,  M.  Sanchez in \cite{sanchez}) in reference to the Lorentzian metric with the general relativity background.  
\end{remark}


\section{Perturbing by shear vector field}

Analyses involving Randers spaces are generally difficult and finding solutions to the geodesic equations is not straightforward \cite{brody, chern_shen}. That is why by solving the problem under the particular perturbation, i.e. shear vector field, we show the respective steps of the solution to obtain at the end the flows of Randers geodesics representing the time-optimal paths. Next, we present  the final solutions graphically which are obtained by following analogous steps for other river-type perturbations - the Gaussian function and the quartic curve perturbation.   

The general solution refers to the system of geodesic equations $\eqref{geo1}$ applied to the Randers metric. We create some computational programmes by means of Wolfram Mathematica ver. 10.2 to generate the graphs and evaluate some numerical computations if the complete symbolic ones cannot be obtained. Using our software in the case of Finslerian computations can strengthen our intuition as well as ckeck and simplify some of obtained formulae and expressions. 
\subsection{The resulting metric}
\begin{figure}
        \centering
~\includegraphics[width=0.33\textwidth]{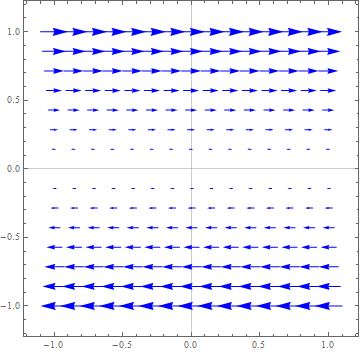}
~\includegraphics[width=0.40\textwidth]{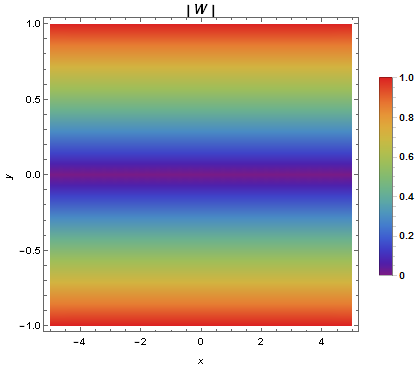}
        \caption{Shear vector field and its density plot under convexity restriction, $|y|<1$.} 
\label{pole_fig}
\end{figure}

We perturb $h$ by the shear vector field (Figure $\ref{pole_fig}$) which is also applied to the discrete problems in the optimal control theory. In the example the current of the river increases as a linear function of $y$ reaching its minimal absolute value in the midstream. E. Zermelo desribed the perturbation  \textit{"as the simplest nontrivial example of our theory"} in his variational paper when the navigation problem had been initially formulated. Let
\begin{equation}
	W = (y, 0)  \qquad \text{with} \qquad  |W| = |y|<1. 
\label{pole}
\end{equation}
\noindent
The condition $|y|<1$ ensures that $F$ is strongly convex and it is a necessary codition if we want to study  the problem via the construction of the Randers metric. Let us observe that E. Zermelo \cite{zermelo}, followed by C. Caratheodory \cite{caratheodory}, considered the problem on the open Euclidean sea without such a restriction what makes worth mentioning difference in both approaches. As a consequence is the fact we are forced to omit formally the solutions under stronger current ($ |W| \geq 1$) in Finslerian approach even in the low-dimensional Euclidean landscape.   
\noindent
We obtain the resulting Randers metric for the shear perturbation
\begin{equation}
F(x,y; u, v) = \frac{ \sqrt{(1-y^2) (u^2+v^2) +(uy)^2} } {1-y^2} - \frac{uy} {1-y^2}.
\end{equation}
The metric can also be presented in the form $F = \alpha + \beta$ as the sum of the new Riemannian metric and 1-form, i.e. 
\begin{itemize}
\item the new Riemannian metric
\begin{equation}
\alpha(x,\mathbf{y})  = \frac{\sqrt {(\delta_{ij}\lambda+W^iW^j)y^iy^j}}{\lambda}
=\frac{\sqrt {(\lambda+(W^1)^2)(y^1)^2+(\lambda+(W^2)^2)(y^2)^2+2W^1W^2y^1y^2}}{\lambda},
\end{equation}
\noindent
and adopting 2D notations yields
\begin{equation}
\alpha(x,y; u, v)= \frac{\sqrt{u^2+v^2-(uW^2-vW^1)^2}}{1-|W|^2}=\frac{ \sqrt{ (1-y^2)v^2 +u^2 } } {1- y^2},
\label{alpha}
\end{equation}
\item 1-form
\begin{equation}
\beta(x,y; u, v)=\frac{-1} {1-|W|^2}(W^1u+W^2v) = \frac{-yu}{ 1 - y^2 }.
\end{equation}
\end{itemize}

\noindent
Hence, the Randers metric in the problem is expressed in the final form by 
\begin{equation}
\label{moja}
F(x,y; u, v) = \alpha(x,y; u, v) + \beta(x,y; u, v) =\frac{ \sqrt{u^2 + v^2 -( yv)^2}-yu } {1- y^2}. 
\end{equation}
Under the influence of $W$, the most efficient navigational paths are no longer the geodesics of the Riemannian metric $h:=(\delta_{ij})$. Instead, they are the geodesics of the Finsler metric $F$. 
Z. Shen \cite{chern_shen} showed that the same phenomenon as in $\mathbb{R}^2$ holds for arbitrary Riemannian backgrounds in all dimensions. The following theorem \cite{colleen_shen} lets us formally to combine the Randers geodesics with the optimality condition  
\begin{theorem}{\text{[D. Bao, C. Robles, Z. Shen, 2004]}}
\label{THM}
A strongly convex Finsler metric $F$ is of Randers type if and only of it solves the Zermelo navigation problem on some Riemannian manifolds $(M, h)$, under the influence of a wind $W$ with $h(W,W)<1$. Also, $F$ is Riemannian if and only if $W=0$. 
\end{theorem}
\noindent
Recall that Randers metrics may be identified with solutions to the navigation problem on Riemannnian manifolds. This navigation structure establishes a bijection between Randers spaces $(M, F=\alpha+\beta)$ and pairs $(h, W)$ of Riemannian metrics $h$ and vector fields $W$ on the manifold $M$. Note, basing on Theorem \ref{THM} and the idea on analysing the indicatrices, the solution to the navigation problem has also  been developed on Hermitian manifolds in complex Finsler geometry with application of complex Randers metric \cite{aldea}.   


\subsection{Projective flatness}

A Finsler metric $F=F(x,\mathbf{y})$ on an open subset $\mathcal{U}\subset \mathbb{R}^n$ is said to be projectively flat if all geodesics are straight lines in $\mathcal{U}$. Our background metric is projectively flat as this is the Euclidean metric. We verify the fact that the resulting Finsler metric \eqref{moja} is not projectively flat. Due to the action of the perturbing vector field which we applied the obtained Randers metric looses this property. In general, there are Riemannian metrics $h$ composing the navigation data $(h, W)$ in the Zermelo navigation problem which are projectively flat, for instance the Klein metric. Depending on the type of perturbation the resulting Randers metric may have or may not have this geometric property. A Finsler metric $F=F(x,\mathbf{y})$ on an open subset $\mathcal{U}\subset \mathbb{R}^n$ is projectively flat iff it satisfies the following system of equations \cite{chern_shen}
\begin{equation}
\label{flat}
F_{x^ky^l}y^k-F_{x^l}=0.
\end{equation}
We prove that the metric is not projectively flat. From (\ref{flat}) we observe in dimension two that
\begin{equation}
u F_{xu}+v F_{yu}=\frac{\left(1-y^2\right) \left(\frac{u v^2 y}{\left(u^2-v^2 \left(y^2-1\right)\right)^{3/2}}-1\right)+2 y \left(\frac{u}{\sqrt{u^2-v^2 \left(y^2-1\right)}}-y\right)}{\left(y^2-1\right)^2}v\neq F_x=0 
\end{equation}
and
\begin{equation}
\small{\quad u F_{xv}+v F_{yv}=\frac{v^4 y}{\left(u^2-v^2 \left(y^2-1\right)\right)^{3/2}}\neq F_y=\frac{-u \left(y^2+1\right) \sqrt{u^2-v^2 \left(y^2-1\right)}+2 u^2 y-v^2 y \left(y^2-1\right)}{\left(y^2-1\right)^2 \sqrt{u^2-v^2 \left(y^2-1\right)}}}.
\end{equation}
Thus, the metric $\eqref{moja}$ is not projectively flat Finsler metric due to the fact that $\eqref{flat}$ does not hold. As a consequnce we cannot obtain the spray coefficients as well as the flag curvature by simplified formulae with use of the scalar function, namely  the projective factor $P$ of $F$, where $G^i=Py^i$. In this case we shall compute both quantities applying the formulae for general Finsler metric. Recall also that a Finsler metric $F$ on manifold $M$ is said to be locally projectively flat if at any point, there is a local coordinate system $(x^i)$ in which $F$ is projectively flat. 

We also verify the projective flatness of the new Riemannian metric $\eqref{alpha}$ of the Randers metric $\eqref{moja}$. We show that 
\begin{equation}
u F_{xu}+v F_{yu}=\frac{u v y \left(2 u^2-3 v^2 \left(y^2-1\right)\right)}{\left(y^2-1\right)^2 \left(u^2-v^2 \left(y^2-1\right)\right)^{3/2}}\neq F_x=0 
\end{equation}
and
\begin{equation}
\quad u F_{xv}+v F_{yv}=\frac{v^4 y}{\left(u^2-v^2 \left(y^2-1\right)\right)^{3/2}}\neq F_y=\frac{y \left(2 u^2-v^2 \left(y^2-1\right)\right)}{\left(y^2-1\right)^2 \sqrt{u^2-v^2 \left(y^2-1\right)}}.
\end{equation}
\noindent
Similarly as the resulting Randers metric the condition $\eqref{flat}$ is not fulfilled, so the Riemannian component $\alpha$ is not projectively flat, either. The example of a pair $(F, \alpha)$, where both metrics, i.e. the resulting Finslerian and the corresponding new Riemannian, are projectively flat, is the Funk and the Klein metric, respectively. 
\subsection{Flag curvature}

For a tangent plane $\tilde{\pi}\subset T_xM$ containing $\mathbf{y}$, $\mathbf{K}=\mathbf{K}(\tilde{\pi},\mathbf{y})$ is the flag curvature, $\tilde{\pi}=span \{\mathbf{y}, \tilde{u}\}$, where $\tilde{u}\in \tilde{\pi}$. In dimension two $\tilde{\pi}=T_xM$ is the tangent plane. Thus the flag curvature $\mathbf{K}=\mathbf{K}(x,y)$ is a scalar function on $TM\backslash\{0\}$. For a Finsler metric $F=F(x,y;u,v)$ the Gauss curvature $\mathbf{K}=\mathbf{K}(x,y;u,v)$ is given by 
\begin{equation}
\mathbf{K}=\frac{1}{F^2}(-2 G_v H_u+2 G Q_u-G_u^2+2 G_x+2 H Q_v-H_v^2+2 H_y-u Q_x-v Q_y)
\end{equation}
where $Q=G_u+H_v$ and $G=G(x,y;u,v)$, $H=H(x,y;u,v)$ denote its spray coefficients. The formula is a formula for the Ricci scalar $Ric$ divided by $F^2$. In 2D the quotient $Ric/F^2$ is the Gauss curvature. We obtain the formula \eqref{Gauss_curvature} which expresses the curvature of the Randers metric in the navigation problem with the shear perturbation
\begin{equation}
\label{Gauss_curvature}
\mathbf{K}=-3\frac{k_n}{k_d}
\end{equation}
where 
\begin{equation}
\scriptsize{
\begin{split}
k_d=4 \sqrt{u^2-v^2 \left(y^2-1\right)}\left(\sqrt{u^2-v^2 \left(y^2-1\right)}-u y\right)^2 \\ 
\left(y \left(y^2+3\right) u^3-\left(3 y^2+1\right) \sqrt{u^2-v^2 \left(y^2-1\right)} u^2-3 v^2 y \left(y^2-1\right) u+v^2 \left(y^2-1\right) \sqrt{u^2-v^2 \left(y^2-1\right)}\right)^4
\end{split}
}
\end{equation}
and
\begin{equation}
\scriptsize{
\begin{split}
k_n=-4 y \left(7 y^{12}+182 y^{10}+1001 y^8+1716 y^6+1001 y^4+182 y^2+7\right) u^{15}\\
+2 \left(y^{14}+91 y^{12}+1001 y^{10}+3003 y^8+3003 y^6+1001 y^4+91 y^2+1\right) \sqrt{u^2-v^2 \left(y^2-1\right)} u^{14}\\
+32 v^2 y \left(2 y^{14}+63 y^{12}+364 y^{10}+429 y^8-286 y^6-455 y^4-112 y^2-5\right) u^{13}\\
+v^2 \left(-3 y^{16}-374 y^{14}-4914 y^{12}-14014 y^{10}-3432 y^8+14014 y^6+7826 y^4+886 y^2+11\right) \sqrt{u^2-v^2 \left(y^2-1\right)} u^{12}\\
-2 v^4 y \left(y^2-1\right)^2 \left(23 y^{12}+1056 y^{10}+8921 y^8+21648 y^6+16929 y^4+3968 y^2+191\right) u^{11}\\
+v^4 \left(y^2-1\right)^2 \left(y^{14}+241 y^{12}+4895 y^{10}+23199 y^8+33495 y^6+15191 y^4+1801 y^2+25\right) \sqrt{u^2-v^2 \left(y^2-1\right)} u^{10}\\
+2 v^6 y \left(y^2-1\right)^3 \left(5 y^{12}+440 y^{10}+5379 y^8+16896 y^6+16115 y^4+4440 y^2+245\right) u^9\\
-3 v^6 \left(y^2-1\right)^3 \left(15 y^{12}+605 y^{10}+4224 y^8+7986 y^6+4455 y^4+625 y^2+10\right) \sqrt{u^2-v^2 \left(y^2-1\right)} u^8\\
-4 v^8 y \left(y^2-1\right)^4 \left(30 y^{10}+723 y^8+3330 y^6+4133 y^4+1390 y^2+90\right) u^7\\
+2 v^8 \left(y^2-1\right)^4 \left(105 y^{10}+1491 y^8+4197 y^6+3075 y^4+530 y^2+10\right) \sqrt{u^2-v^2 \left(y^2-1\right)} u^6\\
+4 v^{10} y \left(y^2-1\right)^5 \left(63 y^8+588 y^6+1080 y^4+472 y^2+37\right) u^5\\
-v^{10} \left(y^2-1\right)^5 \left(210 y^8+1245 y^6+1367 y^4+307 y^2+7\right) \sqrt{u^2-v^2 \left(y^2-1\right)} u^4\\
-2 v^{12} y \left(y^2-1\right)^6 \left(60 y^6+235 y^4+152 y^2+15\right) u^3\\
+v^{12} \left(y^2-1\right)^6 \left(45 y^6+113 y^4+37 y^2+1\right) \sqrt{u^2-v^2 \left(y^2-1\right)} u^2\\
+2 v^{14} y \left(y^2-1\right)^7 \left(5 y^4+8 y^2+1\right) u-v^{14} y^2 \left(y^2-1\right)^7 \left(y^2+1\right) \sqrt{u^2-v^2 \left(y^2-1\right)}.
\end{split}
}
\end{equation}

\begin{figure}
        \centering
~\includegraphics[width=0.4\textwidth]{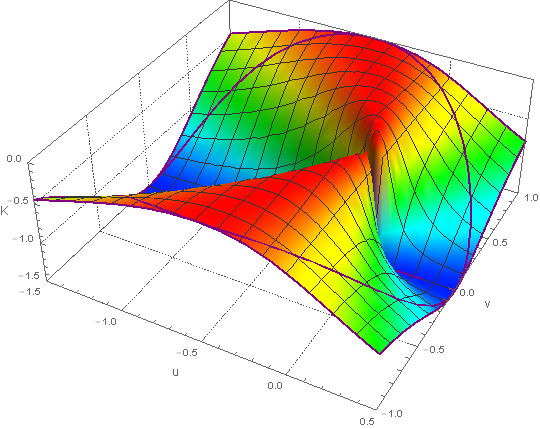}
~\includegraphics[width=0.05\textwidth]{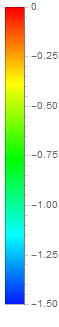}
\includegraphics[width=0.3\textwidth]{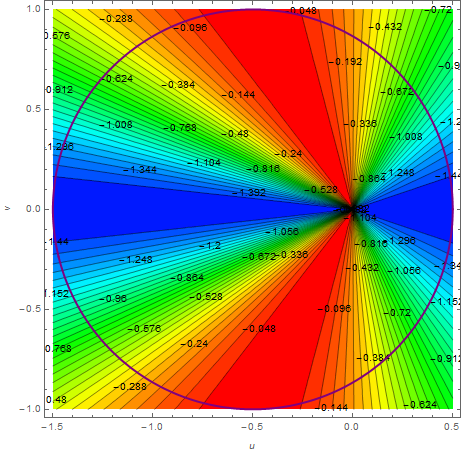}
        \caption{The flag (Gaussian) curvature $\mathbf{K}$ at $A=(0, -\frac{1}{2})$ and its contour plot.}
\label{fig-K}
\end{figure}
\begin{figure}
        \centering
~\includegraphics[width=0.35\textwidth]{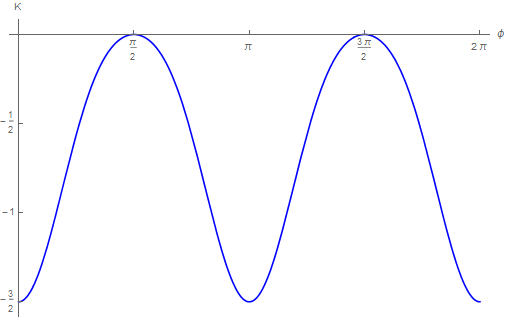}
        \caption{The flag (Gaussian) curvature $\mathbf{K}$ as a function of the initial control angle $\varphi_0\in[0, 2\pi) $ of rotating tangent vector, $A=(0, -\frac{1}{2})$.}
\label{fig-K2}
\end{figure}

Let us look at some computed values of the curvature of the Randers metric. For instance, we obtain $\mathbf{K}(0,-\frac{1}{2};0,\frac{\sqrt{3}}{2})=\mathbf{K}(0, 0; \pm \frac{1}{2},\pm\frac{\sqrt{3}}{2})= -\frac{15}{64}\approx-0.234$; $\mathbf{K}(0, 0; \frac{1}{\sqrt{2}},\frac{1}{\sqrt{2}}) = -\frac{9}{16}\approx -0.563 $. 
In general, the range is determined by $x\in \mathbb{R} \  \wedge \ y\in (-1, 1)$. Taking into consideration the type of the perturbation the ranges refer to the point $(x, y)$ which we choose and $v\in [-1, 1], u\in (y_0-1, y_0+1) \Rightarrow u\in (-2, 2)$. The flag curvature does not depend on $x$-coordinate due to form of the current parallel to the banks of the river. Let us fix $A=(x_0, y_0)=(0, -\frac{1}{2})$ which represents the starting point of the time-optimal paths. The graph of $\mathbf{K}$ at  $A$ for corresponding $u\in (-\frac{3}{2}, \frac{1}{2})$ and $v\in [-1, 1]$ as well as its contour plot are  presented in Figure \ref{fig-K}. However the domain of the tangent vector's coordinates in the problem is restricted and it is represented graphically by the space curve (purple) lying on the 3D surface of the curvature. In other words the length of the resulting tangent vector to the Randers geodesic at the commence point is determined by the equations of motion while the ship's own speed is unit. Gaussian curvature changing given as the function of the initial angle of rotating tangent vector in the tangent space is shown in Figure \ref{fig-K2}. 
In fixed position the minimal value of $\mathbf{K}$ equals $-\frac{{3}}{2}$ and at maximum it approaches 0. The Randers metric \eqref{moja} is not of constant flag curvature. 

\subsection{In steady flow of a current}

\label{constant}
For steady current we choose any two constants $p, q$ which satisfy $p^2 + q^2 < 1$. Let us assume for a while that the flow of the river is represented by the particular river-type vector field $W=[p, 0]$, $|W|=|p|<1$ due to convexity of $F$. Then the resulting Randers metric $\eqref{W1W2}$ is expressed by 
\begin{equation}
F(x,y; u,v) = \frac{ \sqrt{u^2 + v^2 -( pv)^2}-pu } {1-p^2}. 
\end{equation}
This metric is of Minkowski type. The geodesics are the straight lines what can be intuitevely expected. The least time traverse trajectories have a constant velocity and control, i.e. heading angle $\varphi$. The flag curvature is constant, $\mathbf{K}=0$ and the metric fulfills the conditions for projective flatness with the projective factor $P=0$. More general, for two dimensional background with constant vector field, $W=[p, q]$ with $ |W|=\sqrt{p^2+q^2}<1$ the geodesics are obviously the straight lines, too. Thus, the ship sails along the time-optimal path following the same straight track and constant course over ground before and after perturbation while the courses through the water (headings) differ due to acting perturbation. As a consequence there are different resulting speeds over ground as well as the times of the motions between given point of departure and point of destination. Obviously, since we assume the constant perturbation the problem is invariant under translation. We draw attention to this case. Although it represents the simple scenario in the navigation problem theoretically, it is widely applied even in non-optimal and non-Finslerian geometric models and computational algorithms in real navigational software, what we bring closer in the further part of the paper. 


\section{Geodesics of Randers metric as the solutions to Zermelo's problem}

\subsection{Fundamental tensor}
The function $F$ is positive on the manifold $TM\setminus 0$ whose points are of the form $(x, \mathbf{y})$ with $0 \neq \mathbf{y} \in T_x M $ . Over each point $(x, \mathbf{y})$ of $TM\setminus 0 $ treated as a parameter space we designate the vector space $T_x M$ as a fiber and name the resulting vector bundle $\pi^ \ast TM$ . There is a canonical symmetric bilinear form $g_{ij} dx^i \otimes dx^j$ on the fibers of $\pi^\ast T M$ with
\begin{equation}
g_{ij} := \frac 1 2 [ F^2 ]_{y^i y^j}=\frac{1}{2}\frac{\partial^2F^2}{\partial y^i \partial y^j}.
\label{g_ij_Randers}
\end{equation}
The subscripts $y^i , y^j$ signify here partial differentiation and the matrix $(g_{ij} )$ is the fundamental tensor. We compute the partial derivatives of the quadrature of \eqref{moja} 
\begin{equation}
F^2_u=\frac{2 \left(\frac{u}{\sqrt{u^2-v^2 y^2+v^2}}-y\right) \left(\sqrt{u^2-v^2 y^2+v^2}-u y\right)}{\left(1-y^2\right)^2},
\end{equation}

\begin{equation}
F^2_v=\frac{\left(2 v-2 v y^2\right) \left(\sqrt{u^2-v^2 y^2+v^2}-u y\right)}{\left(1-y^2\right)^2 \sqrt{u^2-v^2 y^2+v^2}},
\end{equation}

\begin{equation}
\small{F^2_{uu}=\frac{2 \left(-2 u^3 y+u^2 \left(y^2+1\right) \sqrt{u^2-v^2 \left(y^2-1\right)}-v^2 \left(y^4-1\right) \sqrt{u^2-v^2 \left(y^2-1\right)}+3 u v^2 y \left(y^2-1\right)\right)}{\left(y^2-1\right)^2 \left(u^2-v^2 \left(y^2-1\right)\right)^{3/2}}},
\end{equation}

\begin{equation}
F^2_{vv}=\frac{2 u^3 y-2 u^2 \sqrt{u^2-v^2 \left(y^2-1\right)}+2 v^2 \left(y^2-1\right) \sqrt{u^2-v^2 \left(y^2-1\right)}}{\left(y^2-1\right) \left(u^2-v^2 \left(y^2-1\right)\right)^{3/2}},
\end{equation}

\begin{equation}
F^2_{uv}=-\frac{2 v^3 y}{\left(u^2-v^2 \left(y^2-1\right)\right)^{3/2}}=F^2_{vu}.
\end{equation}

\noindent
Hence, 2D Hessian matrix of $F$ is given by 

\begin{equation}
\scriptsize{H=\left(
\begin{array}{cc}
 \frac{-2 y u^3+\left(y^2+1\right) \sqrt{u^2-v^2 \left(y^2-1\right)} u^2+3 v^2 y \left(y^2-1\right) u-v^2 \left(y^4-1\right) \sqrt{u^2-v^2 \left(y^2-1\right)}}{\left(y^2-1\right)^2 \left(u^2-v^2 \left(y^2-1\right)\right)^{3/2}} & -\frac{v^3 y}{\left(u^2-v^2 \left(y^2-1\right)\right)^{3/2}} \\ \\
 -\frac{v^3 y}{\left(u^2-v^2 \left(y^2-1\right)\right)^{3/2}} & \frac{y u^3-\sqrt{u^2-v^2 \left(y^2-1\right)} u^2+v^2 \left(y^2-1\right) \sqrt{u^2-v^2 \left(y^2-1\right)}}{\left(y^2-1\right) \left(u^2-v^2 \left(y^2-1\right)\right)^{3/2}} \\
\end{array}
\right)}.
\end{equation}
The determinant of the fundamental tensor equals

\begin{equation}
\footnotesize{\text{det}(H)=\frac{u^3 y \left(y^2+3\right)-u^2 \left(3 y^2+1\right) \sqrt{u^2-v^2 \left(y^2-1\right)}+v^2 \left(y^2-1\right) \sqrt{u^2-v^2 \left(y^2-1\right)}-3 u v^2 y \left(y^2-1\right)}{\left(y^2-1\right)^3 \left(u^2-v^2 \left(y^2-1\right)\right)^{3/2}}}.
\end{equation}


\subsection{Spray coefficients}
A spray on $M$ is a smooth vector field on $TM_0:=TM\backslash \{0\}$ locally expressed in the following form
\begin{equation}
    \label{spray}
G=y^i\frac{\partial}{\partial x^i}-2G^i \frac{\partial}{\partial y^i},
\end{equation}
where $G^i=G^i(x,\mathbf{y})$ are the local functions on $TM_0$ satisfying $G^i(x, \widehat{\lambda} \mathbf{y})=\widehat{\lambda} ^2G^i(x,\mathbf{y}), \quad \widehat{\lambda} >0$. 
The spray is induced by $F$ and the spray coefficients $G^i$ of $G$ given by 
\begin{equation}
    \label{spray2}
G^i:=\frac{1}{4}g^{il}\{ [F^2]_{x^ky^l}y^k-[F^2]_{}x^l\}=\frac{1}{4}g^{il}\left(2\frac{\partial g_{jl}}{\partial x^k}-\frac{\partial g_{jk}}{\partial x^l}\right)y^jy^k,
\end{equation}
are the spray coefficients of $F$. More general, in comparison to 2D components obtained above, plugging  \eqref{Alpha} and \eqref{Beta} in \eqref{g_ij_Randers} one obtains the general Hessian matrix of Randers metric given by 
\begin{equation}
g_{ij}(\mathbf{y})=\frac{F}{\alpha}\left(a_{ij}-\frac{y_i}{\alpha}\frac{y_j}{\alpha}\right)\left(\frac{y_i}{\alpha}+b_i\right)\left(\frac{y_j}{\alpha}+b_j\right)
\end{equation}
or in the equivalent form 
\begin{equation}
g_{ij}(\mathbf{y})=\frac{F}{\alpha}\left[a_{ij}-\frac{y_i}{\alpha}\frac{y_j}{\alpha}+\frac{\alpha}{F}\left(\frac{y_i}{\alpha}+b_i\right)\left(\frac{y_j}{\alpha}+b_j\right)\right]
\end{equation}
where $y_i=a_{ij}y^j$. The inverse $(g^{ij})=(g_{ij})^{-1}$ can be expressed in the following form
\begin{equation}
g^{ij}(\mathbf{y})=\frac{\alpha}{F}a^{ij}+\left(\frac{\alpha}{F}\right)^2\frac{\beta+\alpha b^2}{F}\frac{y^i}{\alpha}\frac{y^j}{\alpha}-\left(\frac{\alpha}{F}\right)^2\left(b^j\frac{y^i}{\alpha}+b^i\frac{y^j}{\alpha}\right)
\end{equation}
where $b^i=a^{ij}b_j$. The determinant of the matrix $(g_{ij})$ is given by
\begin{equation}
\text{det} (g_{ij})=\left(\frac{F}{\alpha}\right)^{n+1}\text{det}(a_{ij}).
\end{equation}
\noindent
Now we consider the spray coefficients for general 2D Finsler metric including non-Euclidean background. Let $$L(x,y;u,v)\text{:=}\frac{1}{2}F^2(x,y;u,v).$$
Hence,
\begin{equation}
G^1=\frac{L_{vv}(L_{xu}u+L_{yu}v-L_x)-L_{uv}(L_{xv}u+L_{yv}v-L_y)}{2(L_{uu}L_{vv}-L_{uv}L_{uv})},
\end{equation}
\begin{equation}
G^2=\frac{-L_{uv}(L_{xu}u+L_{yu}v-L_x)+L_{uu}(L_{xv}u+L_{yv}v-L_y)}{2(L_{uu}L_{vv}-L_{uv}L_{uv})}.
\end{equation}
The spray coefficients $G^1:=G=G(x,y;u,v)$ and $G^2:=H=H(x,y;u,v)$ can be expressed by more suitable formulae \cite{chern_shen}
\begin{equation}
G=\frac{\left(\frac{\partial ^2L}{\partial v^2} \frac{\partial L}{\partial x}-\frac{\partial L}{\partial y} \frac{\partial ^2L}{\partial u\, \partial v}\right)-\frac{\partial L}{\partial v} \left(\frac{\partial ^2L}{\partial x\, \partial v}-\frac{\partial ^2L}{\partial y\, \partial u}\right)}{2 \left[\frac{\partial ^2L}{\partial u^2} \frac{\partial ^2L}{\partial v^2}-\left(\frac{\partial ^2L}{\partial u\, \partial v}\right)^2\right]},
\end{equation}
\begin{equation}
H=\frac{\left(\frac{\partial ^2L}{\partial u^2} \frac{\partial L}{\partial y}-\frac{\partial L}{\partial x} \frac{\partial ^2L}{\partial u\, \partial v}\right)+\frac{\partial L}{\partial u} \left(\frac{\partial ^2L}{\partial x\, \partial v}-\frac{\partial ^2L}{\partial y\, \partial u}\right)}{2 \left[\frac{\partial ^2L}{\partial u^2} \frac{\partial ^2L}{\partial v^2}-\left(\frac{\partial ^2L}{\partial u\, \partial v}\right)^2\right]}
\end{equation}
which we apply to compute the spray coefficients of considered Randers metric. We obtain

\begin{equation}
\label{G}
\scriptsize{G=-\frac{v \left(\sqrt{u^2-v^2 \left(y^2-1\right)}-u y\right)^2 \left(4 u y \sqrt{u^2-v^2 \left(y^2-1\right)}-2 u^2 \left(y^2+1\right)+v^2 \left(y^4-1\right)\right)}{2 \left(y^2-1\right) \left(u^3 (-y) \left(y^2+3\right)+u^2 \left(3 y^2+1\right) \sqrt{u^2-v^2 \left(y^2-1\right)}-v^2 \left(y^2-1\right) \sqrt{u^2-v^2 \left(y^2-1\right)}+3 u v^2 y \left(y^2-1\right)\right)}},
\end{equation}

\begin{equation}
\label{H}
\tiny{H=\frac{\left(\sqrt{u^2-v^2 \left(y^2-1\right)}-u y\right)^2 \left(u^3 \left(-\left(3 y^2+1\right)\right)+u^2 y \left(y^2+3\right) \sqrt{u^2-v^2 \left(y^2-1\right)}-v^2 y^3 \left(y^2-1\right) \sqrt{u^2-v^2 \left(y^2-1\right)}+3 u v^2 y^2 \left(y^2-1\right)\right)}{2 \left(y^2-1\right)^2 \left(u^3 y \left(y^2+3\right)-u^2 \left(3 y^2+1\right) \sqrt{u^2-v^2 \left(y^2-1\right)}+v^2 \left(y^2-1\right) \sqrt{u^2-v^2 \left(y^2-1\right)}-3 u v^2 y \left(y^2-1\right)\right)}}. 
\end{equation}

\noindent
Additionally, we also obtain the spray coefficients $G_{\alpha}, H_{\alpha}$ of the new Riemannian term \eqref{alpha} of the resulting Randers metric
\begin{equation}
G_{\alpha}= -\frac{2 u v y}{y^2-1}, \qquad 
H_{\alpha}= \frac{y \left(-2 u^2-v^2 \left(y^2-1\right)\right)}{2 \left(y^2-1\right)^2}. 
\end{equation}


\subsection{Randers geodesic equations}

The solution curves in the Zermelo navigation problem are found by working out the geodesics of the Randers metric. For a standard local coordinate system $(x^i, y^i)$ in $TM_0$ the geodesic equation for  Finsler metric is expressed in the general form
\begin{equation}
\label{geo}
\dot{y}^i+2G^i(x,\mathbf{y})=0. 
\end{equation}
Hence, 
\begin{equation}
\label{geo1}
\ddot{x}^i+\frac{1}{2} g^{il} (2\frac{\partial g_{jl}}{\partial x^k} - \frac{\partial g_{jk}}{\partial x^l})\dot{x}^j\dot{x}^k=0.
\end{equation}

\noindent
Let us abbreviate $\Phi = \sqrt{\dot x^2-\delta \dot y^2}, \delta = \left(y^2-1\right), \Delta = \left(y^2+3\right), \psi = \left(3 y^2+1\right), \Psi = \dot y^2 \Phi, \sigma = \dot x^2 \Phi$. Then plugging $\eqref{G}$ and $\eqref{H}$ in $\eqref{geo}$ we obtain the system of geodesic equations in the navigation problem under the shear vector field as follows

\begin{equation}
\ddot x-\frac{2 \left(\dot y\left(\Phi-y \dot x\right)^2 \left(4 y \dot x \Phi-2 \left(y^2+1\right) \dot x^2+\left(y^4-1\right) \dot y^2\right)\right)}{2 \delta \left(\psi \sigma+3 y \delta \dot x \dot y^2-\delta \Psi-y \Delta \dot x^3\right)}=0,
\label{g1}
\end{equation}

\begin{equation}
\ddot y+\frac{2 \left(\left(\Phi-y \dot x\right)^2 \left(\delta y^3 \left(-\dot y^2\right) \Phi+3 \delta y^2 \dot x \dot y^2+\Delta y \sigma-\psi \dot x^3\right)\right)}{2 \delta^2 \left(-\psi \sigma-3 y \delta \dot x \dot y^2+\delta \Psi+y \Delta \dot x^3\right)}=0.
\label{g2}
\end{equation}
\ 

\begin{figure}
        \centering
~\includegraphics[width=0.8\textwidth]{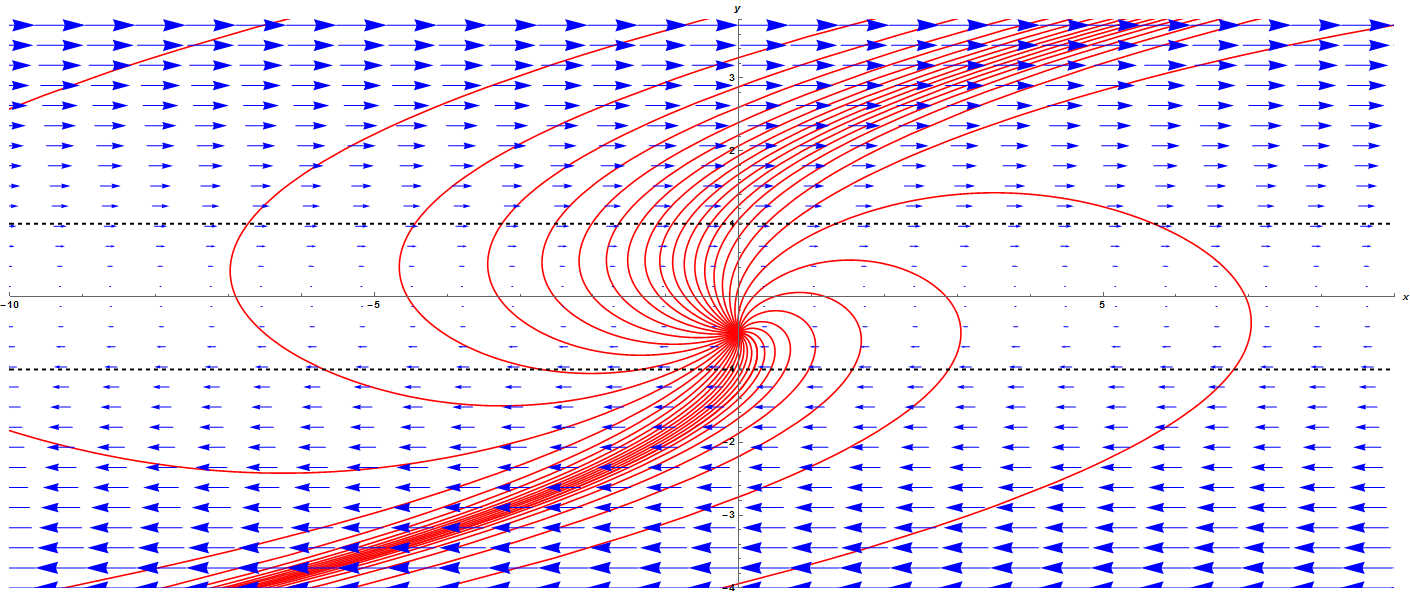}
        \caption{The family of Randers geodesics with the increments $\vartriangle\varphi_0=\frac{\pi}{18}$ at open sea.}
\label{fam_1}
\end{figure}
\begin{figure}
        \centering
~\includegraphics[width=0.57\textwidth]{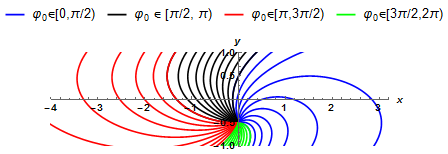}
        \caption{The subsets of Randers geodesics divided with respect to $\frac{\pi}{2}$-ranged $\varphi_0\in[0, 2\pi)$ with the increments  $\vartriangle\varphi_0=\frac{\pi}{18}$ in bounded domain, $|y|<1$.}
\label{fam_1e}
\end{figure}

\noindent
\noindent
Now we evaluate the system of the Randers geodesic equations at full range of the initial angle $\varphi_0$ starting from the fixed point $A=(0, -\frac{1}{2})$. Figure $\ref{fam_1}$ provides graph of the solutions to the system of the equations \eqref{g1} and \eqref{g2} giving the optimal curves for $\varphi_0$ ranging from 0 to 2$\pi$ in the increments of $\vartriangle\varphi_0=\frac{\pi}{18}$ at open sea together with acting vector field. The convexity condition for $|W|$ implies bounded domain, i.e. $|y|<1$, marked by the black dashed lines. 
We also divide the Randers geodesics into the subsets referring to four quarters of $\frac{\pi}{2}$-ranges of the initial angle $\varphi_0$ with $\vartriangle\varphi_0=\frac{\pi}{18}$ in bounded domain what is presented in Figure $\ref{fam_1e}$. The time-optimal paths for $\varphi_0\in[0, \frac{\pi}{2})$ are marked in blue, $\varphi_0\in[\frac{\pi}{2}, \pi)$ in black, $\varphi_0\in[\pi, \frac{3}{2}\pi)$ in red and $\varphi_0\in[\frac{3}{2}\pi, 2\pi)$ in green. The same color codes we shall use in the analogous figures  in the further part of the paper. 
\begin{corollary}{}
\label{cor1}
The obtained results on the navigation problem referring to the mild shear perturbation $\eqref{pole}$ by means of Finsler geometry are consistent to the initial results obtained with the use of Hamiltonian formalism in the calculus of variations (cf. \S 458 in \cite{caratheodory}). 
\end{corollary}
\noindent
In our Finslerian approach the corresponding optimal steering angle can be deduced from obtained equations of Randers geodesics as equivalently from implicit Zermelo's formula applied in $\mathbb{R}^2$.     


\section{Examples of Randers geodesics' flows}

Without loss of generality following above steps of the solution for new Finsler metrics \eqref{W1W2} one  can obtain the final flows of Randers geodesics \eqref{geo} depending on other types of perturbation \eqref{pole_river} and different initial conditions. Next, as acting vector fields we shall consider the ones which are determined by the quartic curve and Gauss function. In Finsler geometry the computations of geometric quantities are usually complicated. As it takes a while if one computes some quotients manually, we create some programmes with use of Wolfram Mathematica to generate the graphs and provide some numeric computations when the complete symbolic ones cannot be obtained. The numerical schemes can give useful information studying the geometric properties of obtained flows. 


\subsection{Quartic curve perturbation}

\begin{figure}
        \centering
~\includegraphics[width=0.35\textwidth]{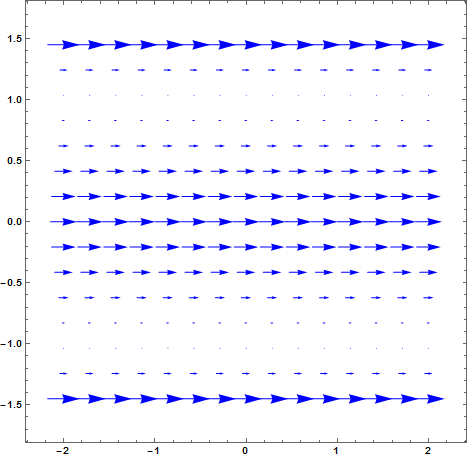}\qquad
~\includegraphics[width=0.15\textwidth]{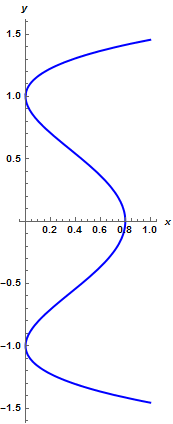}
~\includegraphics[width=0.42\textwidth]{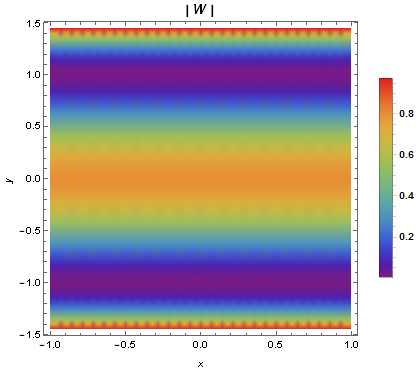}
        \caption{Perturbation determined by the plane quartic curve with $a=0.8$, $b=1$.}
\label{rb_pole}
\end{figure}

\begin{figure}
        \centering
~\includegraphics[width=0.85\textwidth]{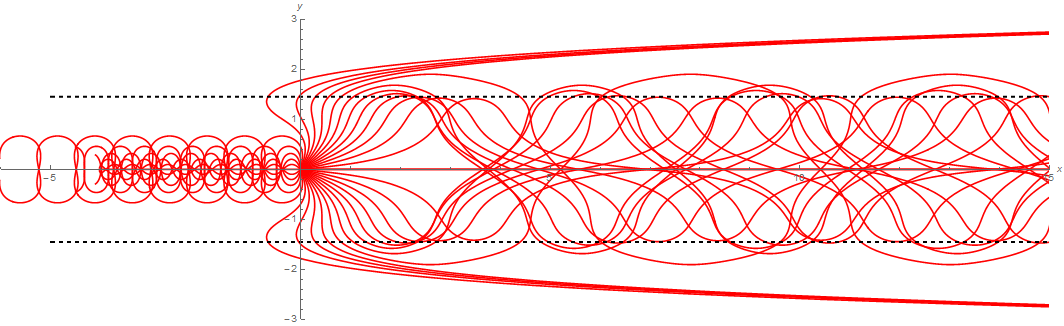}
~\includegraphics[width=0.7\textwidth]{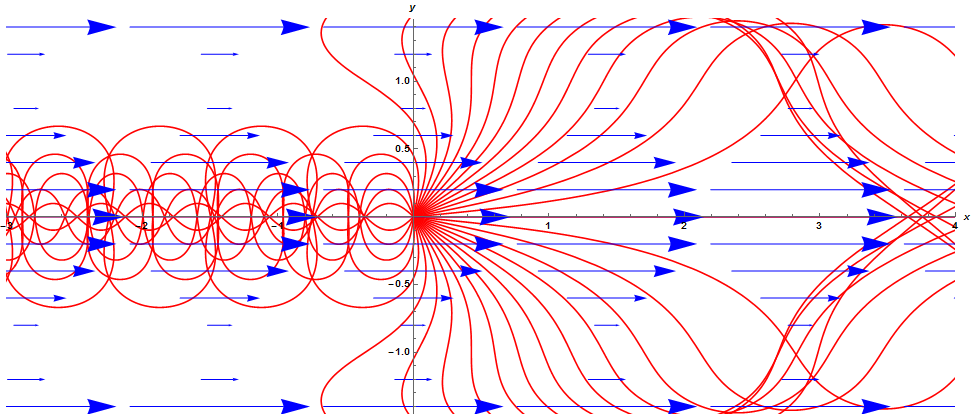}
        \caption{Randers geodesics starting from a point in midstream under quartic curve perturbation in the increments $\bigtriangleup\varphi_0=\frac{\pi}{18}$.}
\label{rb_7c}
\end{figure}

\begin{figure}
        \centering
~\includegraphics[width=0.6\textwidth]{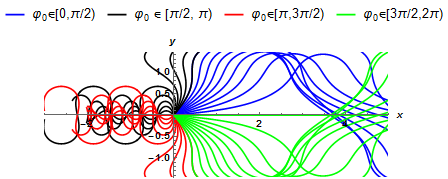}
        \caption{The subsets of Randers geodesics starting from a point in midstream under quartic curve perturbation divided with respect to $\frac{\pi}{2}$-ranged $\varphi_0\in[0, 2\pi)$ in the  increments $\vartriangle\varphi_0=\frac{\pi}{18}$ in bounded domain, $|y|<1.4553$.}
\label{rb_7fffff}
\end{figure}

\begin{figure}
        \centering
~\includegraphics[width=0.9\textwidth]{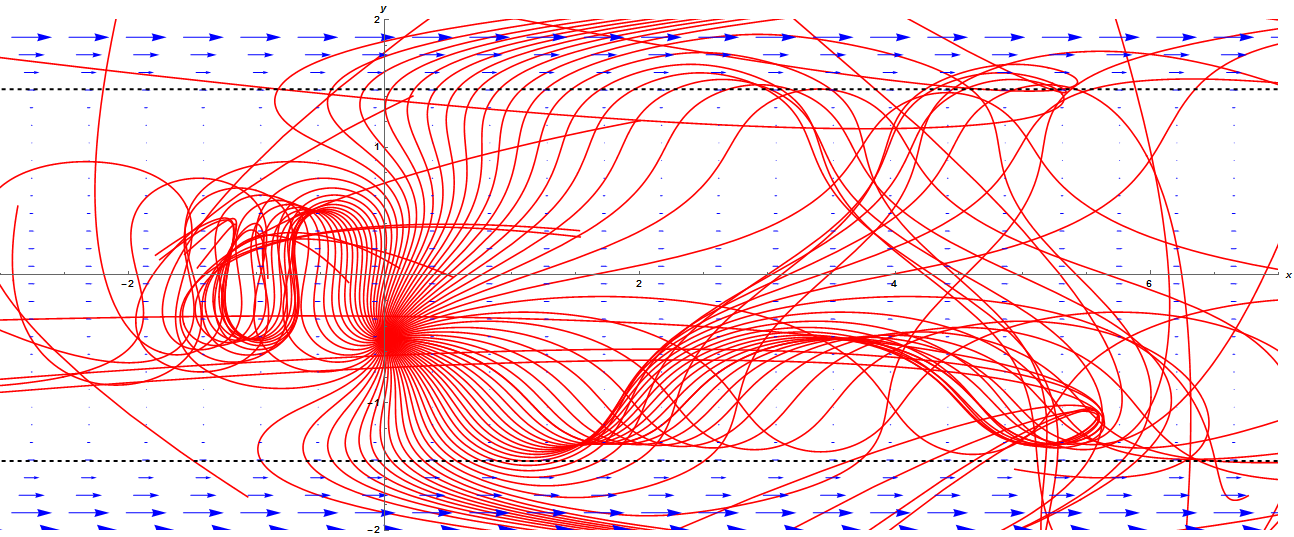}
        \caption{The flow of Randers geodesics starting at $(0, -\frac{1}{2})$ under quartic curve perturbation for $\bigtriangleup\varphi_0=\frac{\pi}{36}$.}
\label{rb_3}
\end{figure}

\begin{figure}
        \centering
~\includegraphics[width=0.6\textwidth]{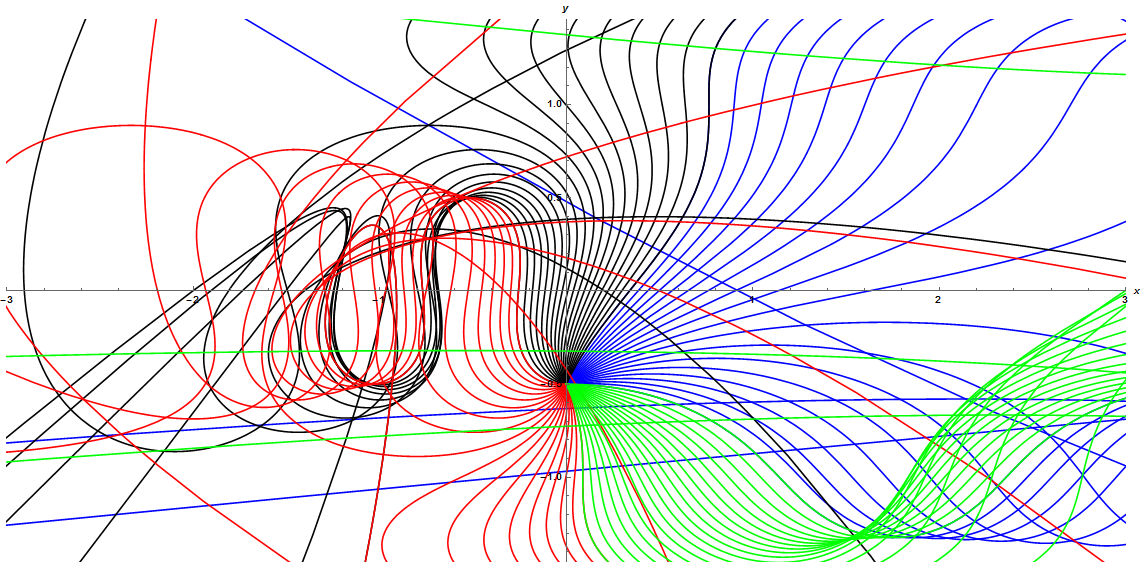}
        \caption{The subsets of Randers geodesics starting at $(0, -\frac{1}{2})$ under quartic curve perturbation divided with respect to $\frac{\pi}{2}$-ranged $\varphi_0\in[0, 2\pi)$  in the increments  $\vartriangle\varphi_0=\frac{\pi}{18}$ in bounded domain, $|y|<1.4553$.}
\label{rb_5}
\end{figure}

\begin{figure}
        \centering
~\includegraphics[width=0.7\textwidth]{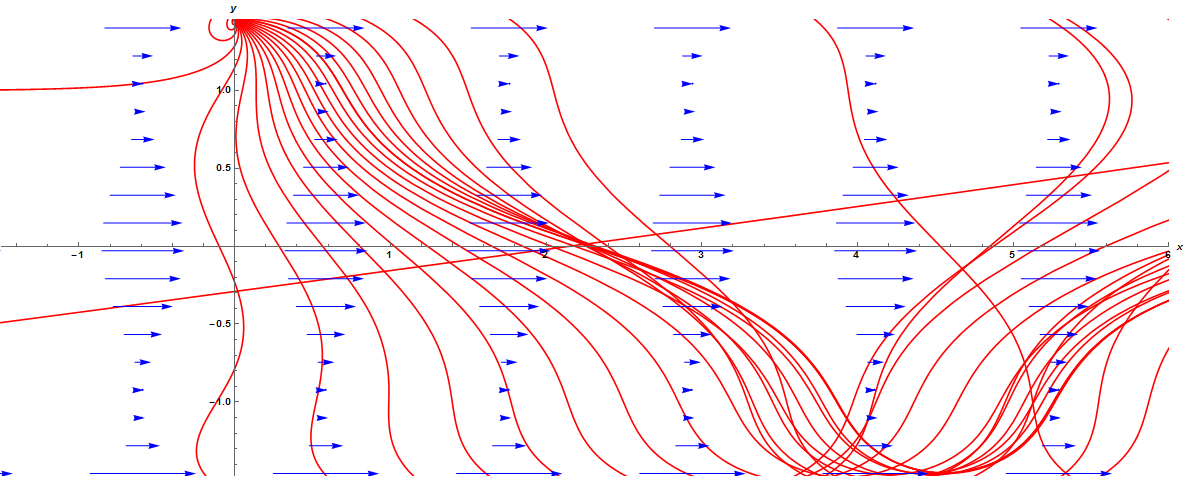}
        \caption{Randers geodesics starting from upper bank of the river at $(0, 1.4553)$ under quartic curve perturbation in bounded domain, $\bigtriangleup\varphi_0=\frac{\pi}{18}$.}
\label{rb_8d}
\end{figure}

\begin{figure}
        \centering
~\includegraphics[width=0.7\textwidth]{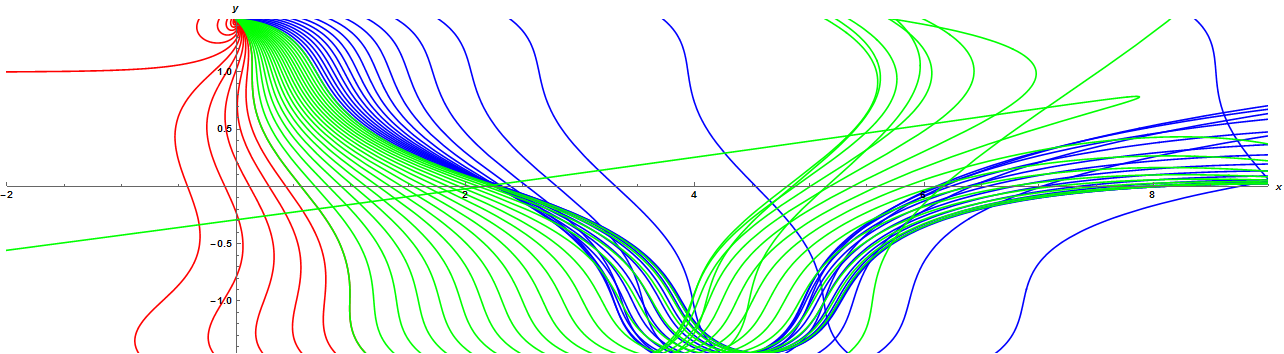}
        \caption{The subsets of Randers geodesics under quartic curve perturbation starting from upper bank of the river at $(0, 1.4553)$, divided with respect to $\frac{\pi}{2}$-ranged $\varphi_0\in[0, 2\pi)$ in bounded domain, with the increments $\vartriangle\varphi_0=\frac{\pi}{18}$.}
\label{rb_8e}
\end{figure}

First let us start with the current $W$ defined by the quartic plane curve given by 
\begin{equation}
f(x)=a \left(b-x^2\right)^2.
\label{quart}
\end{equation}
\begin{figure}
        \centering
~\includegraphics[width=0.41\textwidth]{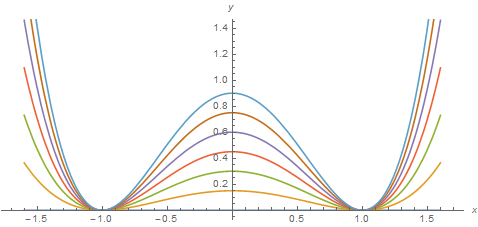}
~\includegraphics[width=0.36\textwidth]{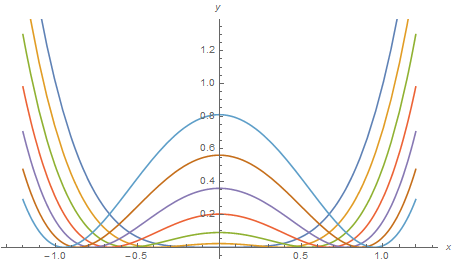}
        \caption{The families of the plane quartic curves \eqref{quart} depending on the parameters $a$, $b$.}
\label{rb_pole2}
\end{figure}
The families of the plane quartic curves for different parameters $a$, $b$ are presented in Figure \ref{rb_pole2}. Let $a=0.8$ and $b=1$ and adopt for the scenario with the horizontal flow of the perturbation. This is shown in Figure \ref{rb_pole}. Recall, to apply Theorem \ref{THM} the convexity condition must be fulfilled. In the set up coordinate system the condition $|f(y)|<1$ implies the domain  $|y|<y_0\approx 1.4553$ which is marked by red dotted lines in Figure \ref{rb_7c}. 
The Randers geodesics starting from a point in midstream generated in the increments $\bigtriangleup\varphi_0=\frac{\pi}{18}$, in the presence of acting quartic curve perturbation are presented in Figure \ref{rb_7c}. Figure \ref{rb_7fffff} shows the subsets of the geodesics divided with  respect to $\frac{\pi}{2}$-ranged $\varphi_0\in[0, 2\pi)$, in the increments  $\vartriangle\varphi_0=\frac{\pi}{18}$, in bounded domain. There are the upstream focal points although there are none downstream. In Figure \ref{rb_3} and Figure \ref{rb_5} we can observe the change of the flow which depends on the initial point from which the family of Randers geodesics leaves. As in previous case of the shear vector field we fix the initial position at $A=(0, -\frac{1}{2})$. The upstream and downstream geodesics connecting the banks of the river and starting from $(0, 1.4553)$ are shown in Figure \ref{rb_8d}, and in partition with respect to $\frac{\pi}{2}$-ranged $\varphi_0\in[0, 2\pi)$ in Figure \ref{rb_8e}.


\subsection{Gaussian function perturbation}

\begin{figure}
        \centering
~\includegraphics[width=0.8\textwidth]{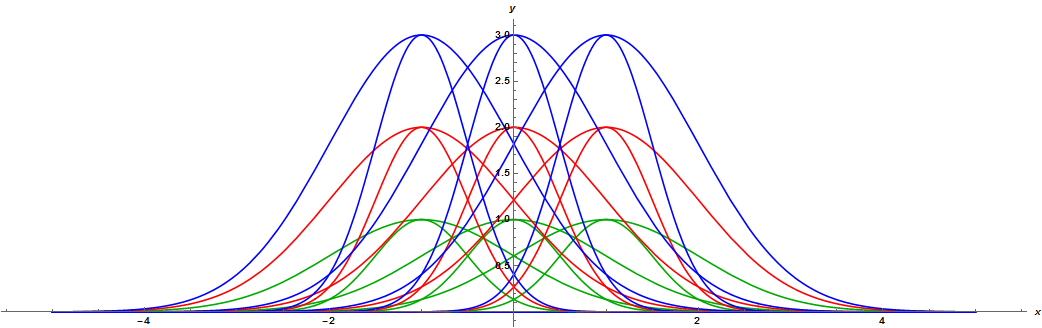}
        \caption{The families of planar Gaussian functions \eqref{gauss} depending on the parameters $a$, $b$, $c$.}
\label{gauss3_5}
\end{figure}

\begin{figure}
        \centering
~\includegraphics[width=0.35\textwidth]{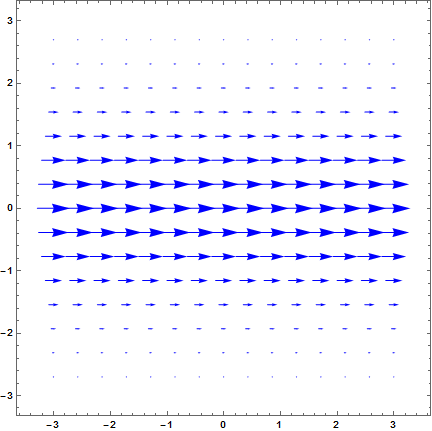}
~\includegraphics[width=0.13\textwidth]{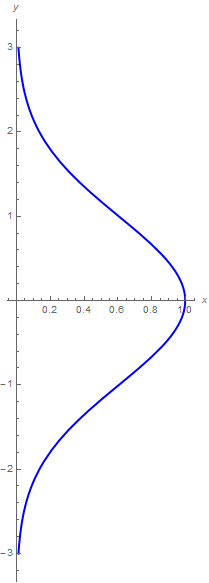}
~\includegraphics[width=0.42\textwidth]{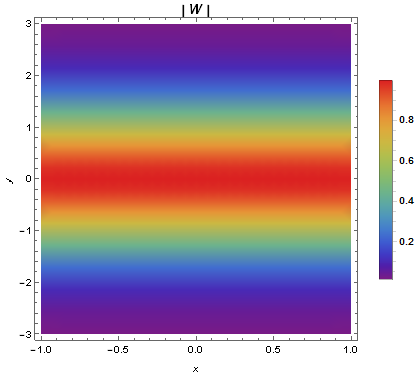}
        \caption{Perturbation determined by Gaussian function with $a=\frac{5}{2\sqrt{2 \pi }}\approx 0.997$, $b=0$, $c=1$.}
\label{gauss3_7}
\end{figure}

\begin{figure}
        \centering
~\includegraphics[width=0.6\textwidth]{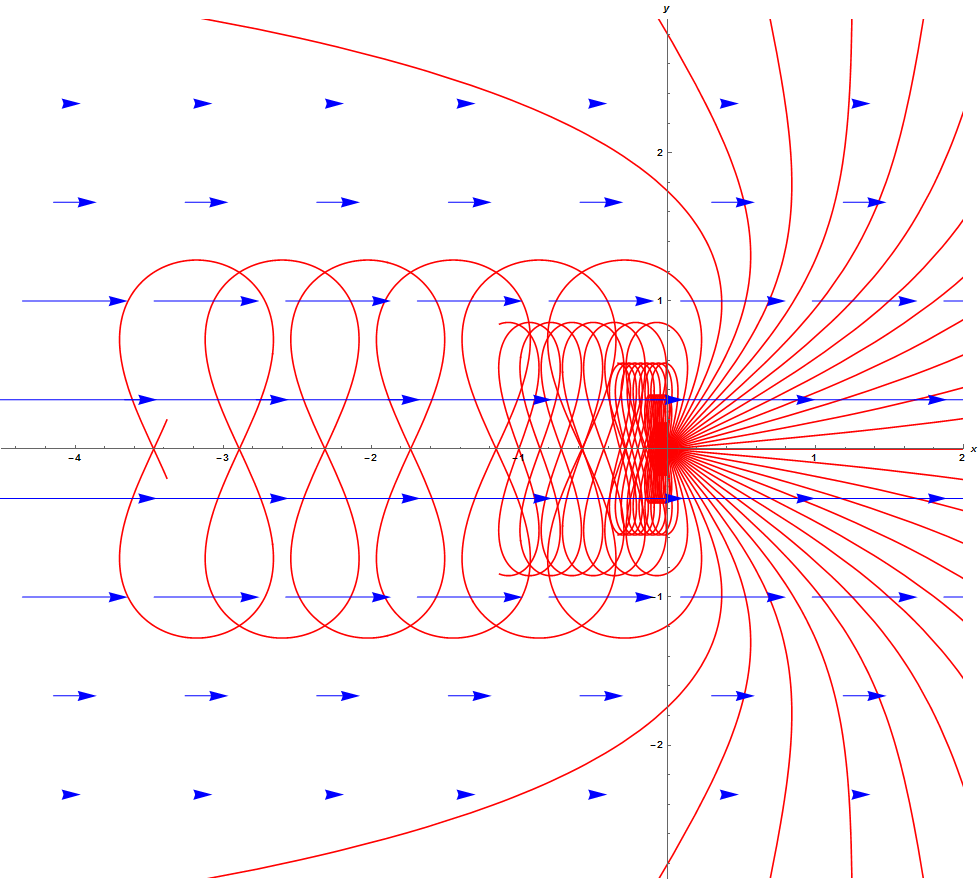}
        \caption{Randers geodesics starting from a point in midstream under Gaussian function perturbation, with the increments $\bigtriangleup\varphi_0=\frac{\pi}{18}$.}
\label{gauss3_10a}
\end{figure}

\begin{figure}
        \centering
~\includegraphics[width=0.6\textwidth]{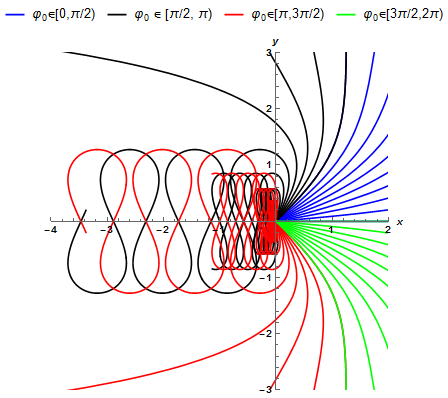}
        \caption{The subfamilies of Randers geodesics starting from a point in midstream under Gaussian function perturbation divided with respect to $\frac{\pi}{2}$-ranged $\varphi_0\in[0, 2\pi)$, with $\vartriangle\varphi_0=\frac{\pi}{18}$.}
\label{gauss3_1b}
\end{figure}

\begin{figure}
        \centering
~\includegraphics[width=0.6\textwidth]{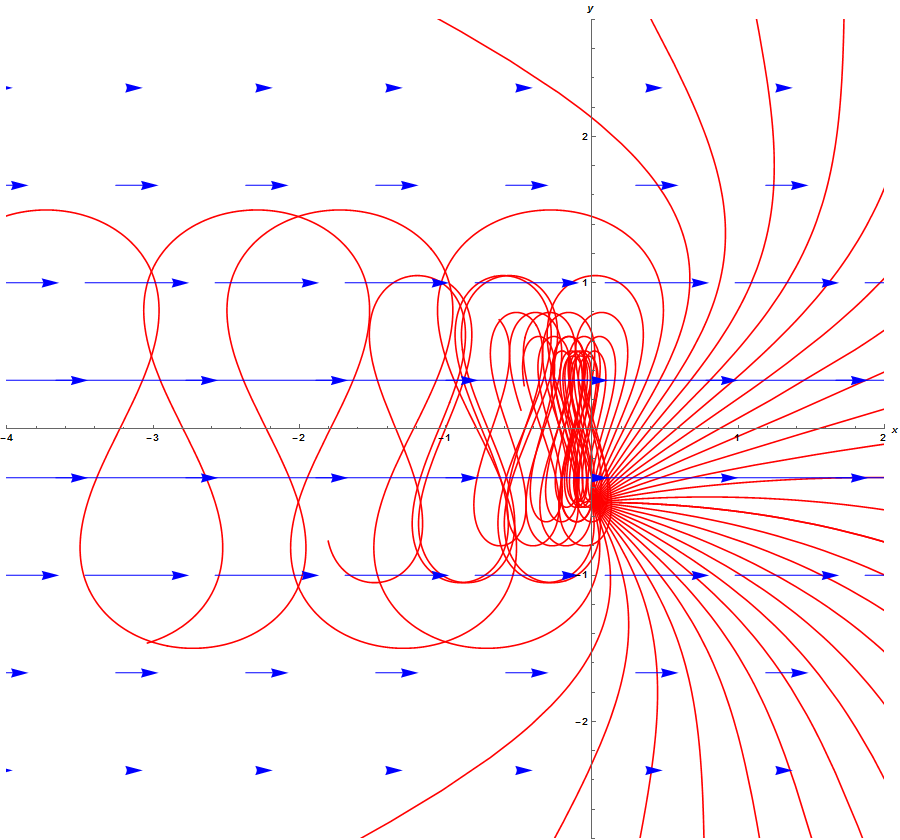}
        \caption{The flow of Randers geodesics starting at $(0, -\frac{1}{2})$ under under Gaussian function perturbation, with the increments $\bigtriangleup\varphi_0=\frac{\pi}{18}$.}
\label{gauss3_9}
\end{figure}

\begin{figure}
        \centering
~\includegraphics[width=0.6\textwidth]{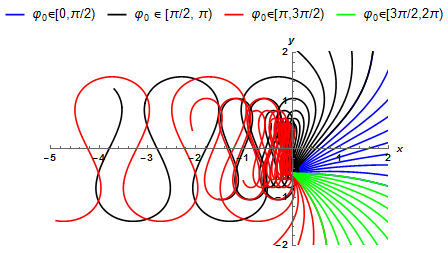}
        \caption{The subsets of Randers geodesics starting at $(0, -\frac{1}{2})$ under Gaussian function perturbation, divided with respect to $\frac{\pi}{2}$-ranged $\varphi_0\in[0, 2\pi)$ in the increments $\vartriangle\varphi_0=\frac{\pi}{18}$.}
\label{gauss3_4AA}
\end{figure}

Next, we consider the current $W$ determined by Gauss function which is expressed in the following form
\begin{equation}
f(x)=ae^{-\frac{(x-b)^2}{2c^2}},
\label{gauss}
\end{equation}
where $a, b, c$ are the arbitrary real constants. The families of plane Gauss function curves for different parameters $a$, $b$, $c$ are presented in Figure \ref{gauss3_5}. This time we want to fullfil the strong convexity condition compulsory for Randers metric applied in the Finslerian approach to Zermelo navigation in the whole space so then the restriction referring to the domain is not necessary. We increase the speed of the current described by the basic Gaussian function $\frac{1}{\sqrt{2 \pi }}e^{-\frac{1}{2}x^2}$ using the constant multiplier, e.g. $\frac{5}{2}$, to approach its maximal allowable value implied by strong convexity. Consequently, let $a=\frac{5}{2\sqrt{2 \pi }}\approx 0.997$, $b=0$, $c=1$ and adopt as before for the scenario with horizontal flow of the perturbation. The field is shown in Figure \ref{gauss3_7}. Hence,
\begin{equation}
f(y)=\frac{5}{2}\frac{1}{\sqrt{2 \pi }}e^{-\frac{1}{2}y^2}.
\end{equation}
\noindent
Time-optimal paths coming from the midstream with the increments  $\bigtriangleup\varphi_0=\frac{\pi}{18}$ in the presence of acting Gaussian function perturbation are presented in Figure \ref{gauss3_10a}. The subfamilies of Randers geodesics starting from a point in midstream and divided with respect to $\frac{\pi}{2}$-ranged $\varphi_0\in[0, 2\pi)$ are presented in Figure \ref{gauss3_1b}, with the increments $\vartriangle\varphi_0=\frac{\pi}{18}$. The analogous flow of Randers geodesics but starting at $(0, -\frac{1}{2})$ with $\bigtriangleup\varphi_0=\frac{\pi}{18}$ is given in Figure \ref{gauss3_9} and Figure \ref{gauss3_4AA}.
\begin{corollary}{}
\label{cor2}
The flows of time-optimal paths in the presence of the mild river-type perturbations $\eqref{pole}$ obtained by means of Finsler geometry are consistent to the results on the Zermelo navigation prolem obtained with the use of the equations of motions and the implicit navigation formula of Zermelo for the optimal control angle in the calculus of variations (cf. \S 279, \S 458 in \cite{caratheodory}). 
\end{corollary}


\section{Application to search models}

\subsection{Revisiting SAR patterns and stating the problem}

Our motivation in current research also touches the search patterns which are applied in air and marine Search and Rescue operations (SAR). They are determined in the International Aeronautical and Maritime Search and Rescue (IAMSAR) Manual which is required to be up-to-date and carried onboard the real ships worldwide by the International Convention for the Safety of Life at Sea (SOLAS, 1974). The Manual provides guidelines for a common aviation and maritime approach to organizing and providing search and rescue services. 
\begin{figure}
        \centering
        \begin{subfigure}{0.8\textwidth}
                \centering
                \includegraphics[width=\textwidth]{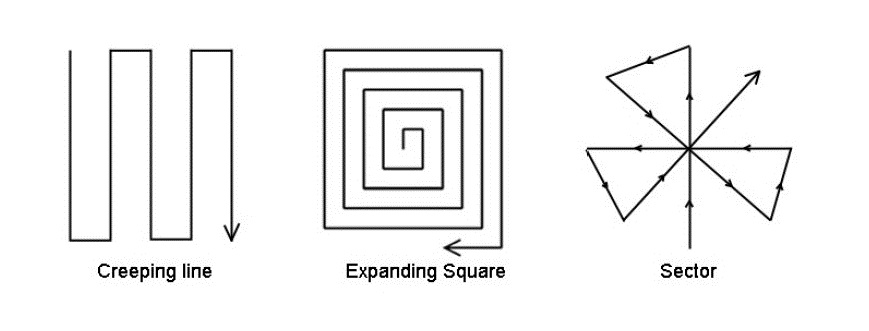}
        \end{subfigure}
   \caption{Planar SAR patterns: creeping line, expanding square and sectror search.}
\label{sar_1}
\end{figure}
\begin{figure}
        \centering
                     \begin{subfigure}{0.4\textwidth}
                \centering
                \includegraphics[width=\textwidth]{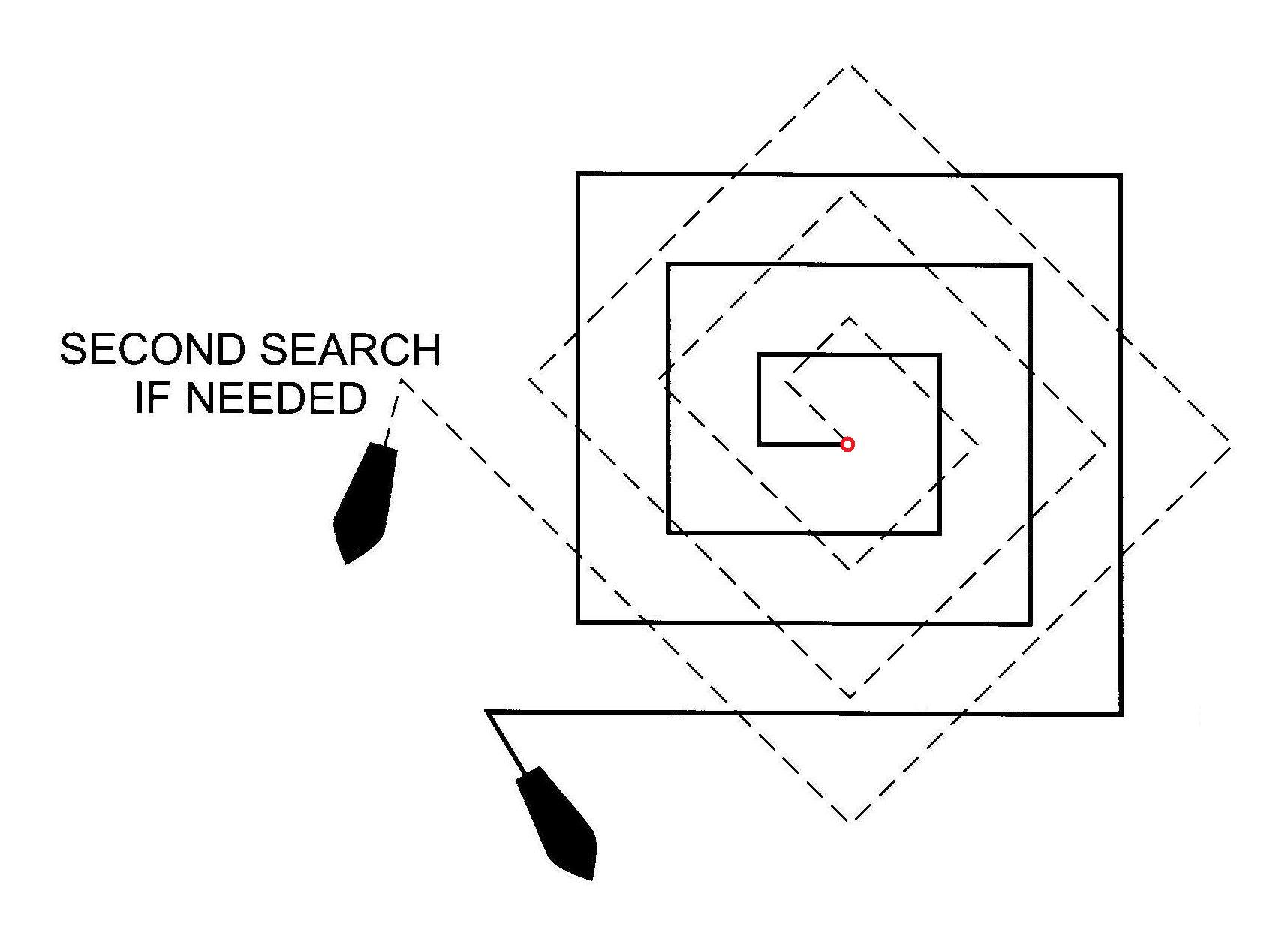}
        \end{subfigure}
              \begin{subfigure}{0.4\textwidth}
                \centering
                \includegraphics[width=\textwidth]{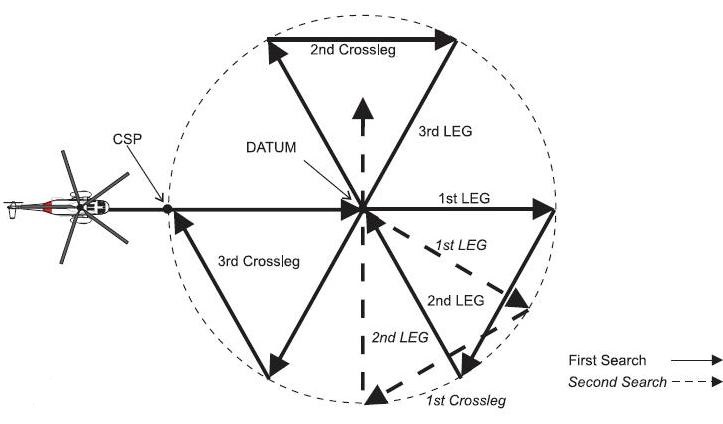}
        \end{subfigure} \quad %

  \begin{subfigure}{0.35\textwidth}
                \centering
                \includegraphics[width=\textwidth]{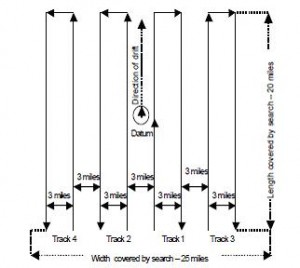}
        \end{subfigure} \quad%
    \begin{subfigure}{0.45\textwidth}
                \centering
                \includegraphics[width=\textwidth]{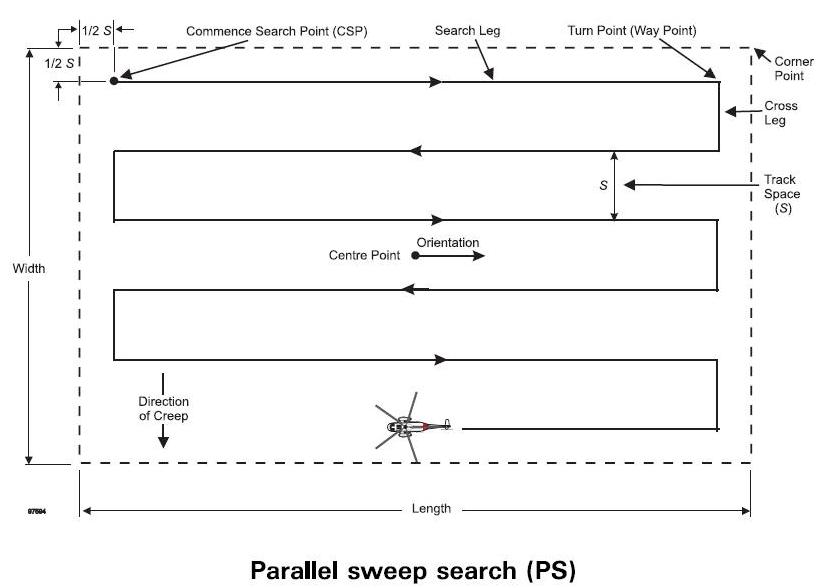}
        \end{subfigure} \quad%
               \caption{Planar SAR patterns in the real air and marine applications: expanding square, sector search, parallel search, creeping line.}
\label{sar_2}
\end{figure}

Briefly, the essential criterium is time in SAR operations. The standard patterns are applied when the accurate position of searched object is unknown. So the key is to search the area efficiently and in the shortest time. Thus, the point of the highest probability (datum), which represents the centre point of searching, is to be determined. In particular, that is the field for applying the theory of reliability and probability including notion of concentration ellipse as in the study of uncertainties and error analysis. However, in our study we focus on the geometry referring to the flows of the time-optimal paths which might be used by searching object. 

We observe that the search models are created and required to be used routinely neglecting the influence of acting perturbation in the meaning of optimization, which in reality may be understood as the real wind or stream (current), in particular. The standard methods of search are based on the following patterns: expanding square, creeping line, sector search and parallel search. They are presented in Figure \ref{sar_1}. Hence, the search paths in the air and marine applications are followed what Figure \ref{sar_2} shows. Now, let us recall two standard geodetic problems in the navigational context from the introduction of the paper. The search path can be sailed in an active or passive way. The former means that we correct the route to follow the search pattern "over ground". In turn, we let the ship to be drifted continuously by acting perturbation in the latter, so we follow the search pattern just "through the water". We ask the question if the search paths, their segments could be based on the time-optimal paths which we study by means of geometry in order to fulfill time criterium more efficient.   


\subsection{Simulations}

In this part of our research we use the navigation information system Navi-Sailor 4000 ECDIS (Electronic Chart Display and Information System) integrated to the navigational simulator Navi-Trainer Professional 5000 by Transas Technologies Ltd., which is one of the worldwide leaders in developing a wide range of IT solutions for the marine industry and promotes its own concept of e-Navigation. We pay our attention to the geometric models and as the consequence the formulae implemented in the navigational software, which the information generated automatically to more and more mass user is based on. In particular, we recall the applications therein which correspond exactly to search models mentioned above. 
\begin{figure}
        \centering
        \begin{subfigure}{0.29\textwidth}
                \centering
                \includegraphics[width=\textwidth]{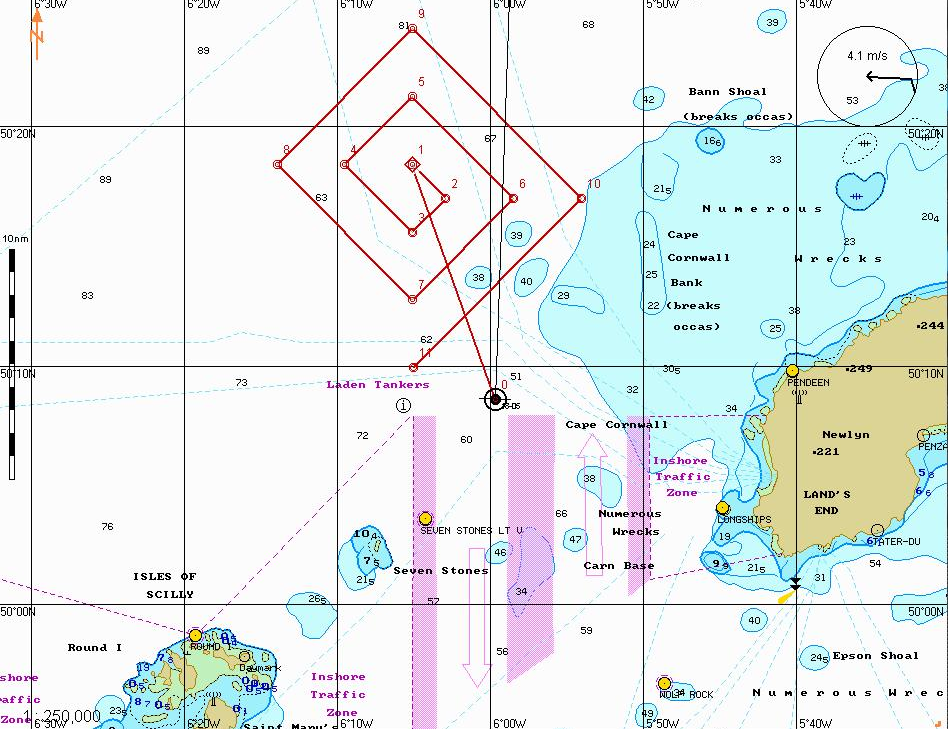}
        \end{subfigure} \quad %
                \begin{subfigure}{0.32\textwidth}
                \centering
                \includegraphics[width=\textwidth]{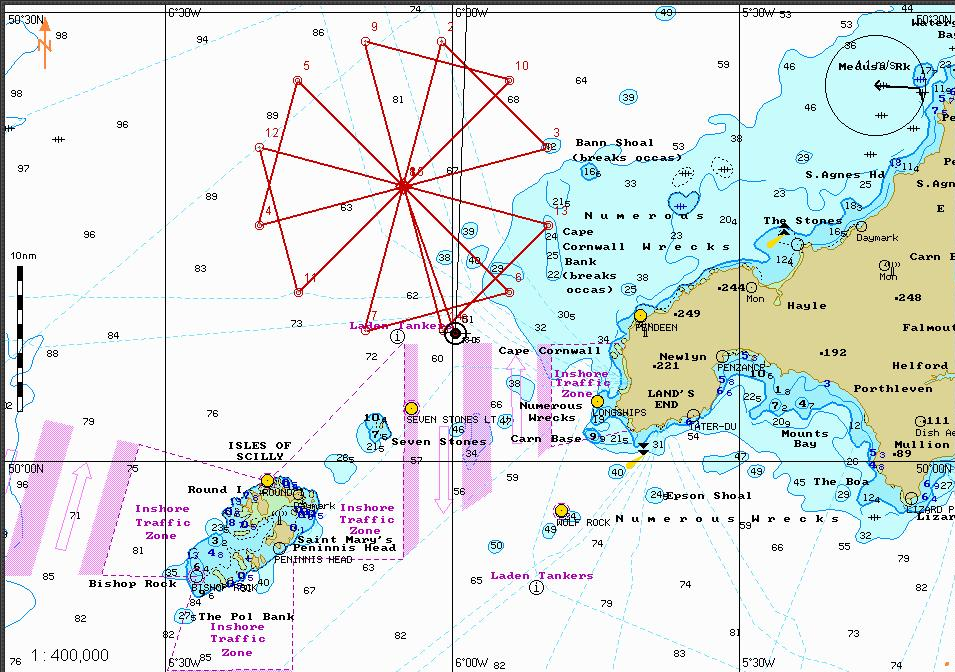}
        \end{subfigure} \quad %
        \begin{subfigure}{0.32\textwidth}
                \centering
                \includegraphics[width=\textwidth]{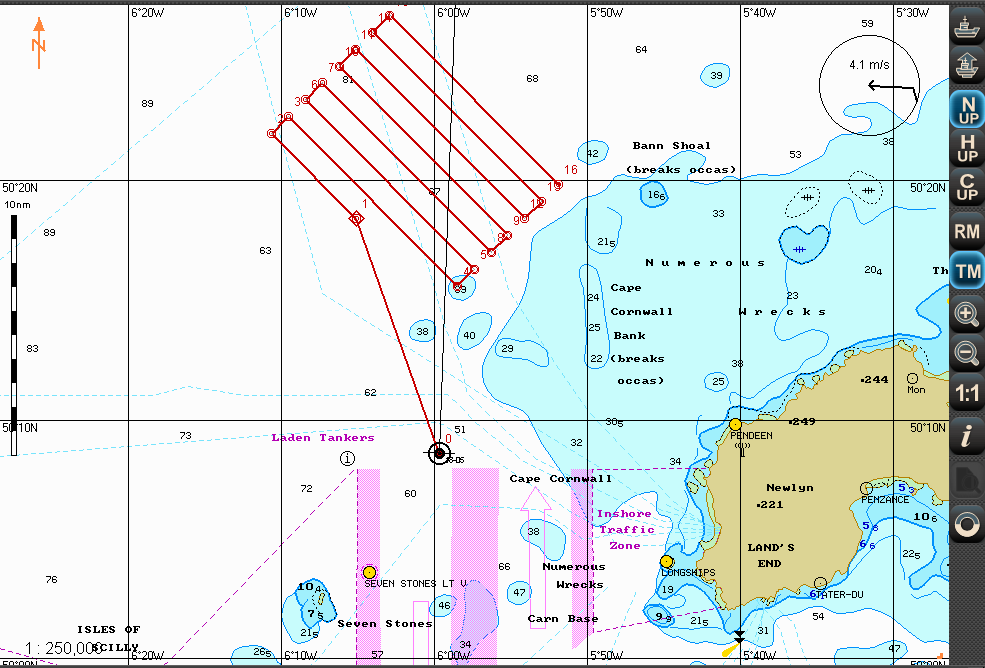}
        \end{subfigure} \quad%
       
        \caption{Simulated search paths in the models of expanding square, sector search and creeping line without perturbation; searching ship positioned in (50$^\circ$08.607'N, 005$^\circ$59.654'W), $|u|=10$ kn.}
\label{sar_bez}
\end{figure}

Let us consider some examples referring to real applications. First, we activate the search paths for the scenario without perturbation (calm sea model). Generated routes based on the models of expanding square, sector search and creepng line in case of one searching ship are presented in Figure \ref{sar_bez}. We assume the same conditions as in the Zermelo navigation problem we considered above, i.e. constant ship's own speed sailing on two-dimensional sea. In practice the search speed means the maximal speed of the ship which is possible to be kept in the real conditions. Note that for $n>1$ searching ships this means the highest speed of the slowest ship. Simulated trajectories correspond to the standard patterns. 
\begin{figure}
        \centering
        \begin{subfigure}{0.32\textwidth}
                \centering
                \includegraphics[width=\textwidth]{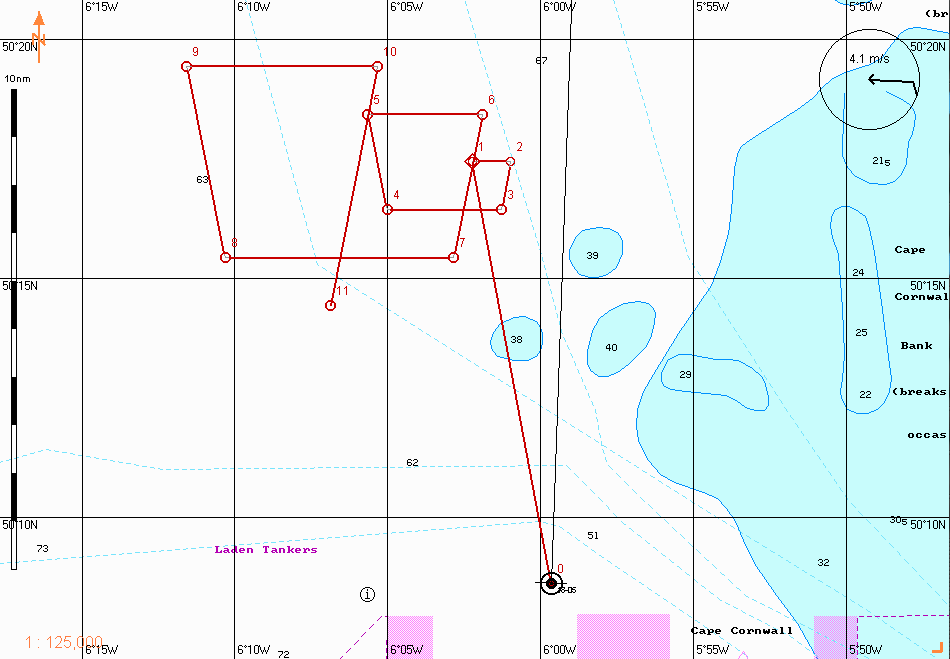}
        \end{subfigure} \quad %
        \begin{subfigure}{0.32\textwidth}
                \centering
                \includegraphics[width=\textwidth]{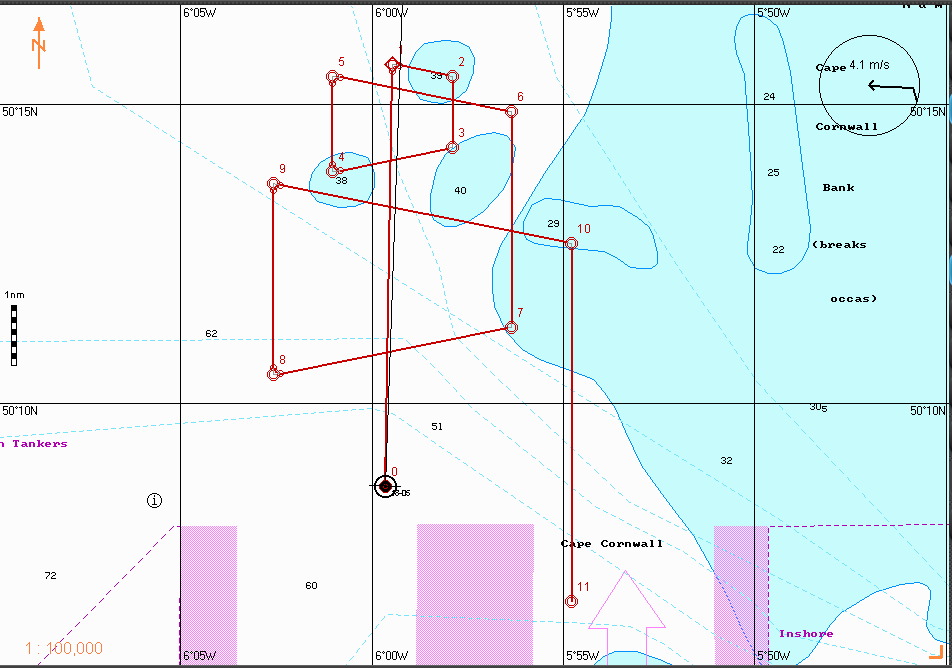}
        \end{subfigure} \quad%
    \begin{subfigure}{0.29\textwidth}
                \centering
                \includegraphics[width=\textwidth]{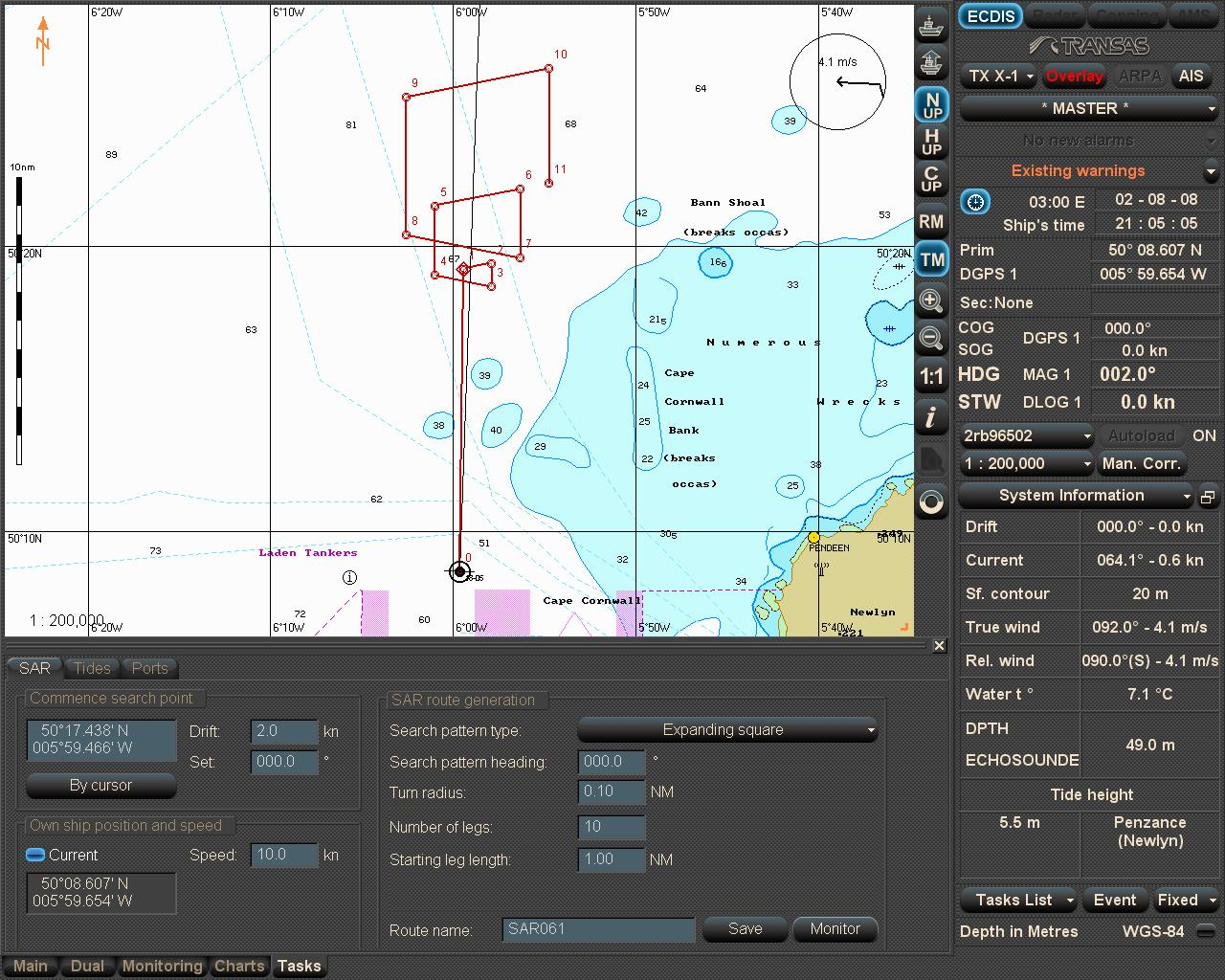}
        \end{subfigure} \quad%
            \caption{Simulated search paths in the models of expanding square under acting perturbations with set of 2.0 kn and drifts: 270.0$^\circ$, 180.0$^\circ$, 000.0$^\circ$, respectively; $|u|=10$ kn.}
\label{sar_z1}
\end{figure}

\begin{figure}
        \centering
        \begin{subfigure}{0.31\textwidth}
                \centering
                \includegraphics[width=\textwidth]{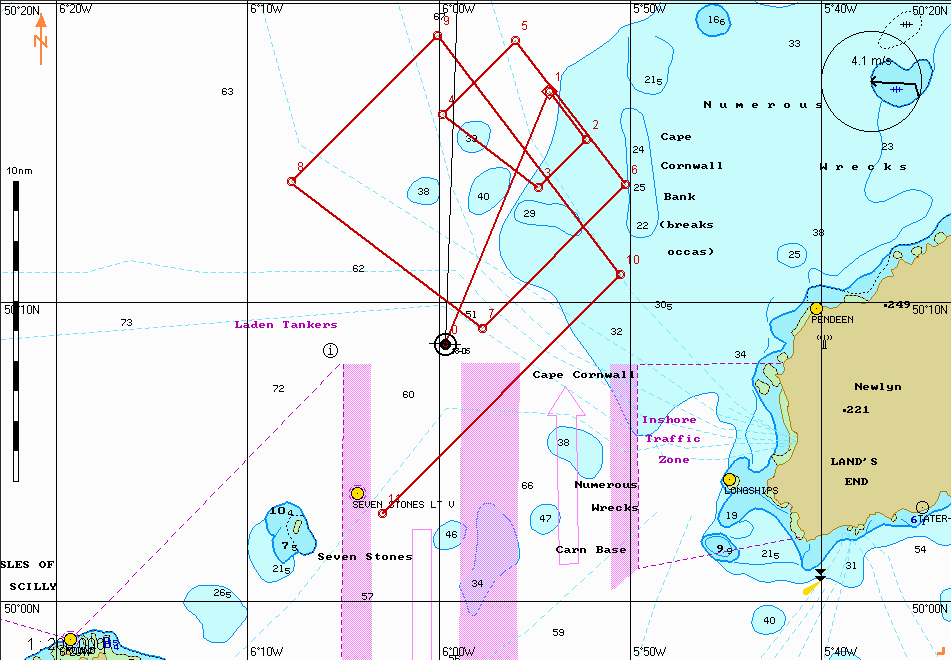}
        \end{subfigure} \quad %
        \begin{subfigure}{0.31\textwidth}
                \centering
                \includegraphics[width=\textwidth]{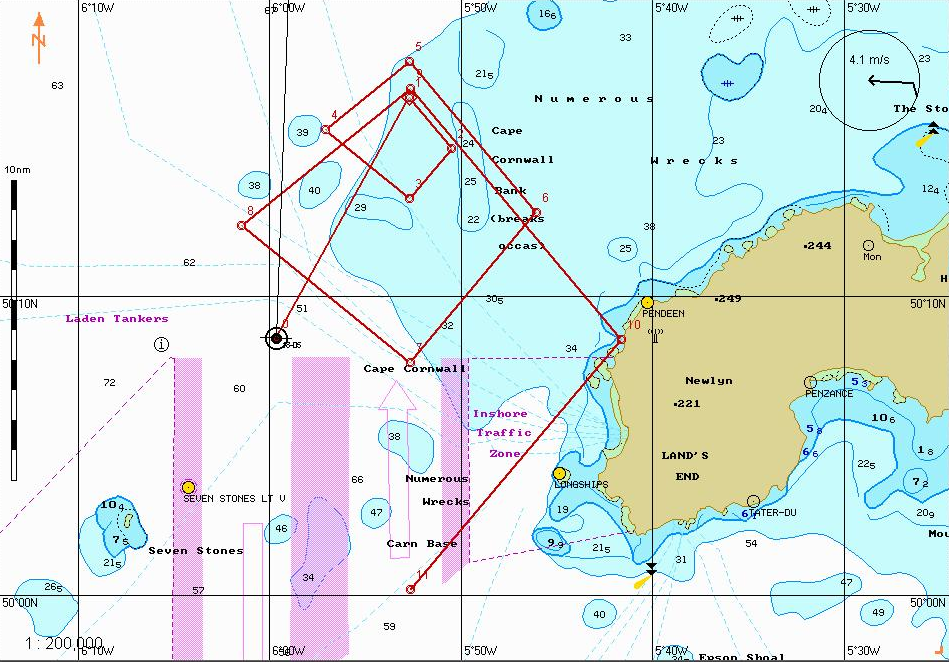}
        \end{subfigure} \quad%
    \begin{subfigure}{0.31\textwidth}
                \centering
                \includegraphics[width=\textwidth]{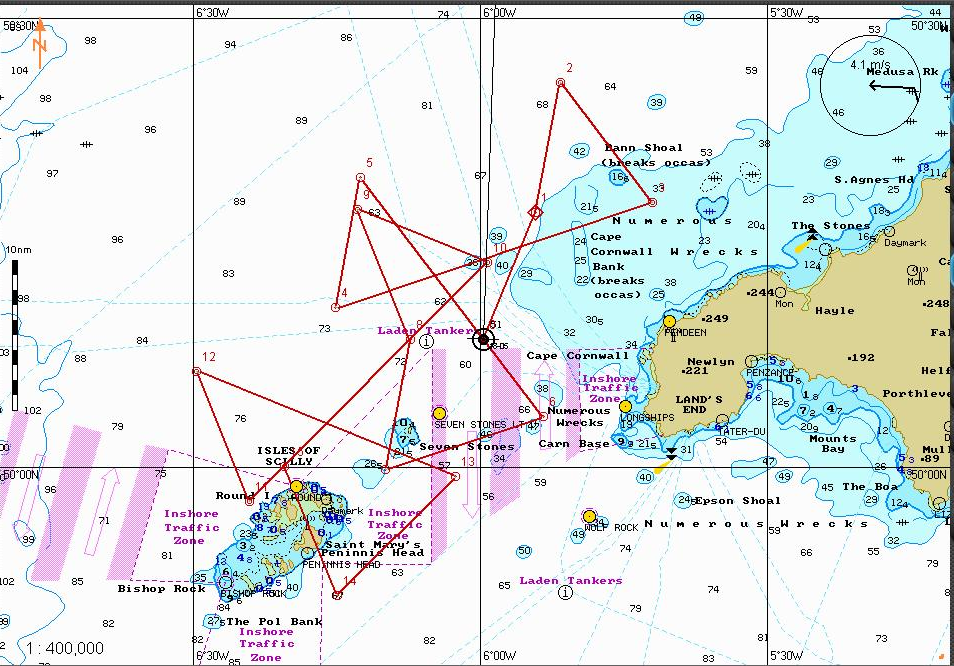}
        \end{subfigure} \quad%
            \caption{Simulated search paths in the models of expanding squares and sector search under "weaker" perturbations with set of 2.0 kn and drifts: 225.0$^\circ$, 180.0$^\circ$, 225.0$^\circ$, respectively; $|u|=10$ kn.}
\label{sar_z2}
\end{figure}

Next, we introduce some perturbations for the same datum in the way it is possible to be set up in the system. The real perturbation can be defined by two parameters, namely the drift (direction) and set (speed).  Let us assume the ship is positioned in fixed point given by the geographical coordinates 50$^\circ$08.607'N (latitude) and 005$^\circ$59.654'W (longitude). We set the following conditions in the simulations. In the first expanding square model presented in Figure \ref{sar_z1} the speed of the ship equals 10 kn (knots), the perturbation is determined by the pairs of drift and set as follows: (270.0$^\circ$, 2.0 kn), (180.0$^\circ$, 2.0 kn), (000.0$^\circ$, 2.0 kn). The commence search point is 50$^\circ$17.438'N and 005$^\circ$59.466'W, number of legs equals 10, the starting leg length is 1 Nm (nautical mile) and the search pattern heading equals 000.0$^\circ$. Then we increase the ship's speed up to 15 kn keeping the same strength of perturbation. Thus, the perturbations become weaker relatively. The solutions obtained in the simulator are given in Figure \ref{sar_z2} and refer to the perturbations given by set and drift as follows (225.0$^\circ$, 2.0 kn), (180.0$^\circ$, 2.0 kn) for expanding square, and (225.0$^\circ$, 2.0 kn) for the sector search. The commence search point is (50$^\circ$18.435'N,  006$^\circ$05.126'W), number of legs equals 10, the search pattern heading is 045.0$^\circ$ and the starting leg length is 2 Nm in case of expanding square. The number of sectors is 6 and the search radius equals 10 Nm in case of the sector search.  

The modified search paths in the models of expanding square and sector search are presented in Figure \ref{sar_z1} and Figure \ref{sar_z2}. Geometrically, the Euclidean plane with steady current is just used in order to generate the search paths. This refers to the case mentioned before in $\S$ \ref{constant}. In fact only constant river-type perturbations are applied in the simulator as well as in the real devices onboard the ships as the same software is used therein. The perturbation is considered in the passive way so it refers to the direct geodetic problem in the navigational context. Due to some practical reasons the simplified approach is followed routinely in reality. However, taking into consideration the state-of-art in positioning, modeling, tracking and control in the context of, e.g. robotic sailing, dron piloting, we may ask if the approach becomes oversimplified as the models can be optimized geometrically. This remark plays a role if the perturbations are variable in space and/or time, so as the real ones.   

\subsection{Preliminary modification with application of time-optimal paths}

\begin{figure}
        \centering
        \begin{subfigure}{0.3\textwidth}
                \centering
                \includegraphics[width=\textwidth]{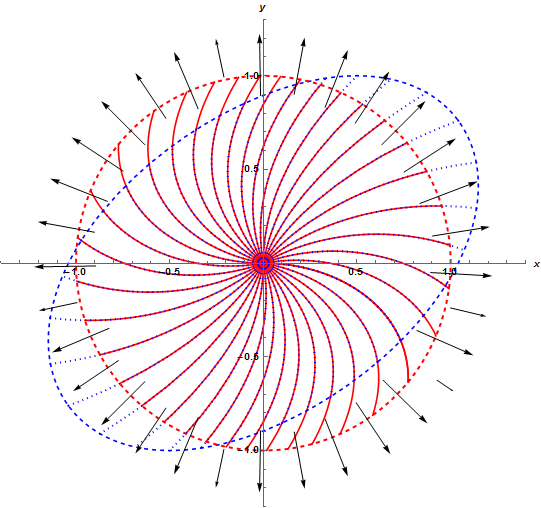}
        \end{subfigure} \quad %
             \begin{subfigure}{0.32\textwidth}
                \centering
                \includegraphics[width=\textwidth]{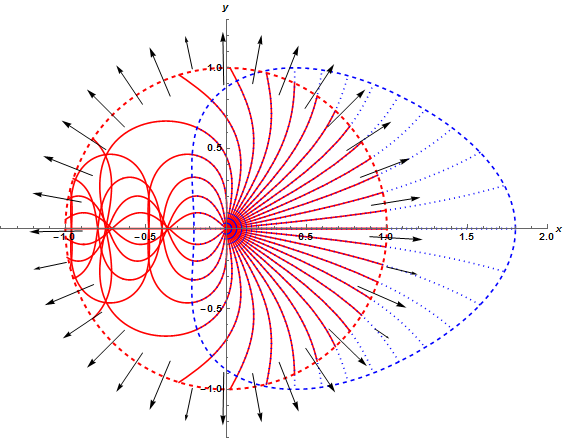}
        \end{subfigure} \quad%
    \begin{subfigure}{0.3\textwidth}
                \centering
                \includegraphics[width=\textwidth]{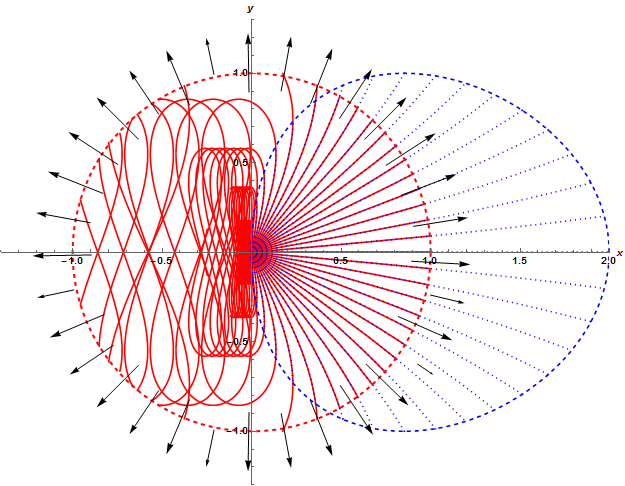}
        \end{subfigure} \quad%
            \caption{Randers geodesics (red solid) in a unit disc starting from the origin, with the increments $\bigtriangleup\varphi_0=\frac{\pi}{18}$ and the indicatrices for $t=1$ (blue dashed) under acting shear,  quartic curve and Gaussian function perturbation, respectively.}
\label{disc}
\end{figure}

Let us first imagine $n$ bees which have been ordered to leave the hive in different directions and reach the boundary of disc-shaped garden in the shortest time such that during the flight in windy conditions the whole hive's neighbourhood is fully patroled. Randers geodesics (red solid) simulating the trajectories in a unit disc starting from the origin, in the increments $\bigtriangleup\varphi_0=\frac{\pi}{18}$, and the corresponding  indicatrices for $t=1$ (blue dashed), under acting the shear, the quartic curve and the Gaussian function background fields, are illustrated in Figure \ref{disc}. In the absence of the winds the optimal trajectories are desribed by the straight rays coming from the origin. They are indicated by the black arrows. In the presence of the wind, for instance in the left-hand side graph of the figure referring to shear perturbation, the space can be fully covered (searched) alternatively by $n$ curved time-optimal paths after suitable adopting the initial angles $\varphi_0$ which determine the flow of the Randers geodesics. The example gives rise to consider the model with only one ship involved in the problem.

In contrast to the search paths created in the software of the navigational simulator let us continue with the model with the perturbation $W$ which is not constant in space. To begin with, we combine the standard patterns with the time-optimal paths which can be represented by the Randers geodesics. Let us remark here that the application of Randers geodesics in the problem can bring more noticeable benefit in solving the problem with non-Euclidean backgrounds where the non-Finslerian approach cannot be applied. Let the current be given by the river-type perturbation \eqref{pole_river}. We consider the piecewise time-optimal paths connecting the fixed waypoints defined in the standard models where the directions of the straight search paths change. With the search problem in mind, we refer to the optimal control angle $\varphi $ by the following corollary 
\begin{corollary}{(cf. \S 459 in \cite{caratheodory})}
The steering must always be toward the side which makes the wind component acting against the steering direction greater. 
\end{corollary}
\noindent
Thus, the idea is to make use of perturbing vector field in order to increase the resulting ship's speed and not to follow the fixed standard pattern without regard for the type and properties of acting perturbation. Brielfy, if possible we aim to avoid sailing routinely "against" the current. 
\begin{figure}
        \centering
~\includegraphics[width=0.145\textwidth]{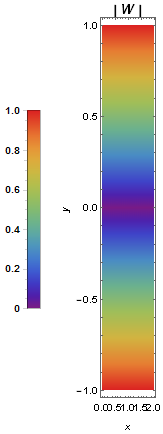} \qquad
~\includegraphics[width=0.5\textwidth]{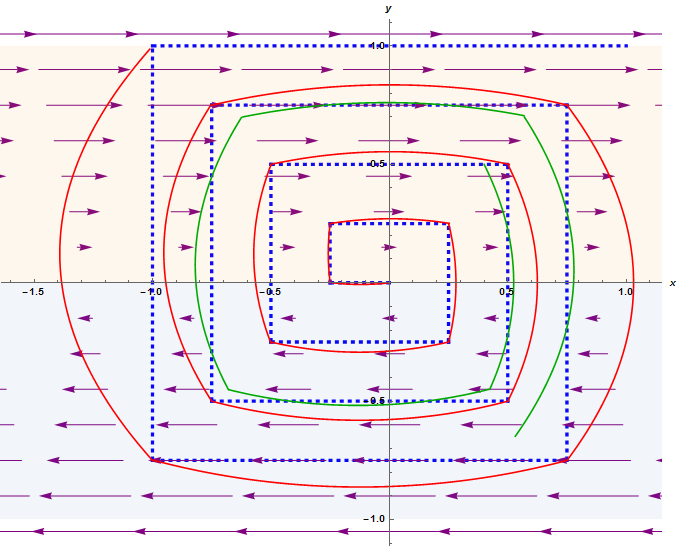}
        \caption{Modified expanding square model by the piecewise-time-optimal legs under the shear vector field.}
\label{exsq}
\end{figure}
The current effect represented by the drift vector is the same for searching and searched object as it is assumed in the simulations. In each of simulated scenarios, if we subtract the drift vector given as linear function of time, then we obtain the standard search patterns illustrated in Figure \ref{sar_bez}. Hence, in reference to flowing water the standard model is followed continuously, however over ground the search paths are then modified by the vector of drift as they are presented in Figure \ref{sar_z1} and Figure \ref{sar_z2}. 

We proceed by first considering expanding square and let us now come back, for simplicity but without loss of general idea, to the weak shear vector field \eqref{pole}. In Figure \ref{exsq} the standard expanding square (blue dotted) of given starting leg length, which determines track spacing $\varepsilon^*$ in the whole model, and oriented such that the horizontal legs are parallel to the flow. Thus, the fixed waypoints are determined and they represent the consecutive startpoints and endpoints connected by the Randers geodesics (red solid). As the time of passage is shorter in each leg in comparison to the corresponding straight legs of the standard pattern, thus the total time in the new modified paths is decreased. Obviously, we require that the search is efficient, so the maximal distance $\varepsilon$ between the points of the searched space and the search path ought to be also taken into consideration. Therefore, in what follows we shall define \textit{the complete search}. So far we aim to show the potential application of (piecewise) time-optimal path, where it is reasonable, in order to minimize the total time of search without exceeding required value of $\varepsilon$. For that reason previously fixed waypoints are now translated such that obtained new Randers geodesics (green solid) fulfill the condition for the maximal distance between the points of searched space and search paths. 
\begin{figure}
        \centering
~\includegraphics[width=0.18\textwidth]{pole2b} \qquad
~\includegraphics[width=0.5\textwidth]{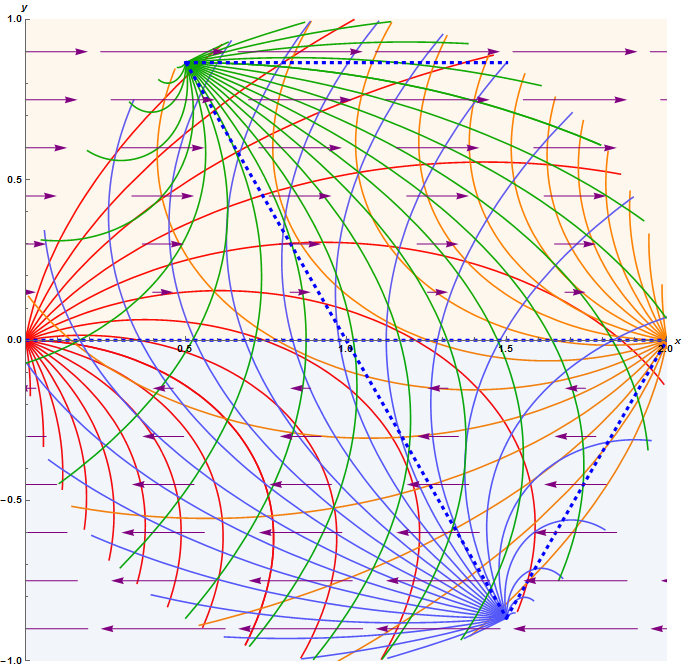}
        \caption{Sector search modified by the time-optimal legs starting from fixed waypoints determined by the standard pattern under the shear vector field, $\bigtriangleup\varphi_0=\frac{\pi}{18}$.}
\label{sector}
\end{figure}

\begin{figure}
        \centering
~\includegraphics[width=0.18\textwidth]{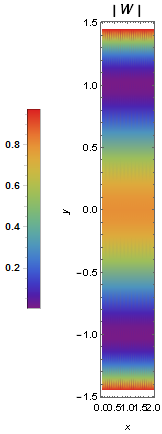} \qquad
~\includegraphics[width=0.45\textwidth]{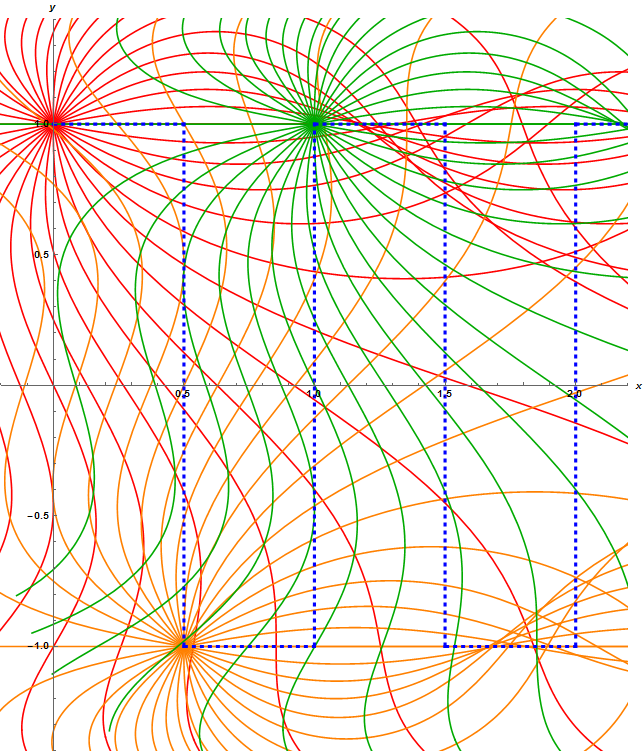}
        \caption{Creeping line search modified by the time-optimal legs starting from fixed waypoints determined by the standard pattern under the quartic curve perturbation,  $\bigtriangleup\varphi_0=\frac{\pi}{18}$.}
\label{creep_RB}
\end{figure}

\begin{figure}
        \centering
~\includegraphics[width=0.18\textwidth]{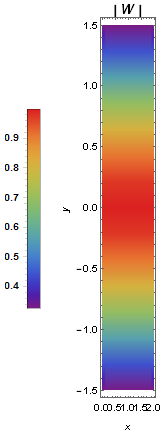} \qquad
~\includegraphics[width=0.5\textwidth]{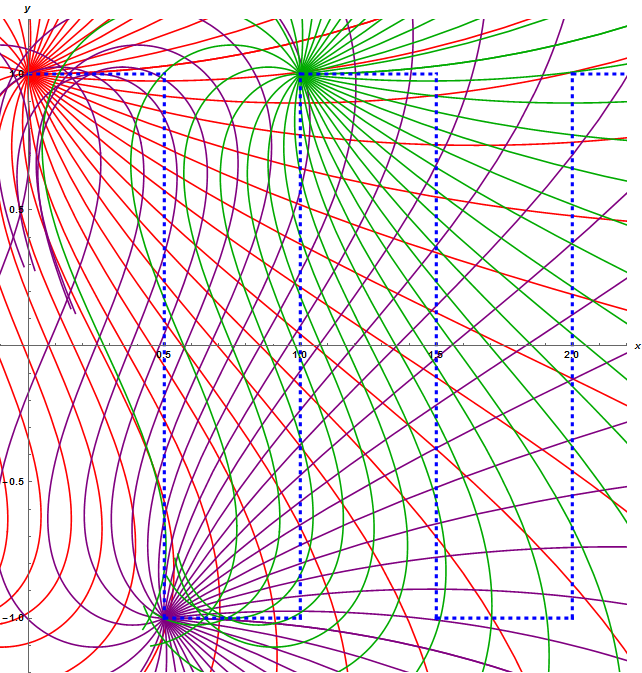}
        \caption{Creeping line modified by the time-optimal legs starting from fixed waypoints determined by the standard pattern under the Gaussian function perturbation, $\bigtriangleup\varphi_0=\frac{\pi}{18}$.}
\label{creep_gauss}
\end{figure}

Next, we combine the Randers geodesics and the sector search under the same weak shear perturbation what is presented in Figure \ref{sector}. The standard pattern (blue dotted) determines, as previously, the consecutive fixed waypoints from which the color-coded families of optimal connections start. In the example we assume that the diameter of circular search area equals the width of the river due to the restriction on strong convexity. If the Randers geodesics connecting directly the fixed points do not fulfill a  condition for $\varepsilon$, then we follow the individuals of the families generated from the intermediate  endpoints of followed legs. In general, the partition depends on $\varepsilon$ and the properties of the time-optimal paths in given vector field. We proceed by considering the third standard pattern, namely the creeping line in the presence of the quartic curve perturbation and the Gaussian function perturbation what is shown in Figure \ref{creep_RB} and Figure \ref{creep_gauss}, respectively. In analogous way as before the standard paths (blue dotted) are modified with the use of the time-optimal legs starting from fixed waypoints determined by the standard pattern in both scenarios. 

The charts of the color-coded Randers geodesics have been created which cover the background spaces under acting river-type perturbations. The curves start from determined fixed points and represent the quickest connections. In this sense they are the solutions to Zermelo's problem in the presence of weak perturbation in entire space without singularities. However, in search problem the additional conditions are  included, e.g. $\varepsilon$ or restricted domain, in the context of free or fixed final time problems which belong to the optimal control theory. In each scenario the conditions modify the final search model. We aim to make use of time-optimal legs to decrease the total time of searching the whole space or to maximize the space which can be searched in limited fixed time, with required value of $\varepsilon$.  
Beside combining the Randers geodesics with the standard search patterns, it can be reasonable in the presence of some perturbations to omit the standard paths completely. This means to follow non-standard search models configured with the time-optimal legs and rearranged wypoints without references  to the standard patterns. The approach admits the perturbation to influence the geometry of the search models. 


\subsection{Generalizations and extensions}
\label{gen}
The problem requires that the search is efficient, not just time-optimal in the meaning of following the time-optimal legs. Thus, introducing the time-optimal paths to the search problem the condition for searching the entire space needs to be fulfilled. This means that the points of the searched space should be close enough to at least one path. In the standard models, which we analyzed above, the track spacing $\varepsilon^*$ is fixed as a constant value and determines the entire model and so the total time of search. Let us define the complete search for a subset $\mathcal{D}$ of a metric space $(M, d)$.  
\begin{definition}
\em{Let $\mathcal{D}\subset M$ be a search space, $\Gamma=\underset{i}{\cup}\gamma_i$ a search path, where $\gamma_i$ represents $i$-time-optimal leg and $\varepsilon\geq0$ is a fixed search parameter. We say that a search is complete if
\begin{equation}
	\forall \ \  A\in \mathcal{D} \ \  \exists \ \ \tilde{A}\in \Gamma: \ d(A, \tilde{A})\leq \varepsilon.
\label{complete_search}
\end{equation}}
\end{definition}
\noindent
The definition implies there are no omitted "zones" left, if the time is free, and the ship follows the time-optimal legs. Hence, the Randers geodesic paths in compliance with the condition \eqref{complete_search} can be applied in the problem. The particular application depends on the initial conditions, type of perturbation and preset parameter $\varepsilon$. In comparison to the track spacing parameter $\varepsilon^*$ applied in the real navigational devices and systems we are led to the relation, for example in case of creeping line search, $\varepsilon=0.5 \varepsilon^*$. However, recall that the key is to minimize the total time $t_c$ of the complete search. Otherwise, the solution based on the piecewise time-optimal search path $\Gamma$, which can be represented by the solutions to Zermelo's problem, might not state for the time-optimal solution to the problem of search. $\gamma_i$ guarantees the local optimality in the connections of the intermediate waypoints. To begin with, we restrict the study excluding the topological singularities, in particular let the search subset $\mathcal{D}$ be compact, connected and convex. The generalized geometric and optimal control problems with $t_c\rightarrow min$, which rise from above, can be proposed as follows
\begin{enumerate}
  \item time-optimal complete search of the ellipse-shaped area under weak river-type perturbation (note, in engineering and navigation the planar positions' distributions of graduated probability levels centered at datum are determined by the error and concentration ellipses, in 3D by the triaxial ellipsoids),  
  \item time-optimal complete search of given 2D area under weak river-type perturbation, arbitrary weak stationary perturbation strong stationary perturbation, arbitrary time-dependent perturbation, 
  \item simultaneous time-optimal complete search of given 2D area by $n\geq2$ ships under above mentioned perturbations,
  \item 3D and higher dimensions analogues of above mentioned tasks with Euclidean background,
  \item analogues of above mentioned tasks with non-Euclidean background, e.g. $\mathbb{S}^n$.
 \end{enumerate}
\noindent
Let us observe that further in the extended problem one can combine and work with the notions of differential geometry and optimization together with reliability and probability. We can ask, for instance, about the complete search model with $\varepsilon\geq0$, in which the searched space is covered by $\Gamma$ such that $t_c\rightarrow min$ while the space is graduated with respect to the probability, starting from the datum (maximum) to the boundary (minimum) of the space. It means that first we search completely the subspace of the highest probability level and then we expand the search with next subspaces of lower probability levels in descending order.   


\section{Conclusions}

We have studied the navigation problem by means of Finsler geometry, namely the Randers metric, under the river-type perturbation. In the case of planar Euclidean background the final system of equations defining the problem can also be formulated and solved with the use of equations of motions and the classical implicit navigation formula of Zermelo for the optimal control angle. Moreover, in order to optimize the real-time computations, what is the essential point in the real applications as the implementation of the geometric models in the software of the real navigational systems and simulators, the latter approach becomes more efficient and simpler. This is due to the fact that in Finsler geometry the computations of geometric quantities are usually complicated, even in 2D with the Euclidean background metric. Analyses involving Randers spaces are generally difficult and finding solutions to the geodesic equations is not straightforward  (cf. \cite{brody, chern_shen}). For instance, this is shown directly in the obtained final system of Randers geodesics' equations \eqref{g1} and \eqref{g2} in the case of the shear wind. In the context of the real applications this is a strong argument, as the computational ability represented by the central processing unit is taken into account in practice, even though each of the mentioned theoretical methods gives the same final solutions represented by the time-optimal paths (cf. Corollary \ref{cor1} \& \ref{cor2}). 

Having chosen the river-type perturbation allows us to simplify the Randers geodesics' equations without loss of generality and present the steps of the solution in completed forms. However, the Randers metric with application of Theorem \ref{THM} is the key to solve the Zermelo navigation problem for non-Euclidean Riemannian background metrics in any dimensions under arbitrary weak perturbations. Thus, the Finslerian approach as well as the variational one \cite{palacek} for an arbitrary wind enable to investigate and solve more advanced geometric and control problems in which, for instance, the implicit formula of Zermelo cannot be used. The analysis of the problem simplifies if we observe that it suffices to find the locally optimal solution on the tangent space. 

The idea behind our scheme is to make use of the time-optimal paths, represented in particular by the Randers geodesics. With the condition \eqref{complete_search} in mind, we aimed to optimize the total time of search in comparison to following the standard fixed search patterns. They are routinely used without taking into consideration the type of the perturbation and the initial conditions. However, these factors influence the geometric properties of the models, the flow of Randers geodesics and thus the resulting speed of searching object, and finally the total traverse time $t_c$. The standard models are also required formally to be used in the real applications. This also motivated us to revisit the standard patterns and focus on the reasons enforcing them to be followed in the scenarios with acting perturbations. In the presented approach we let the the perturbation to influence the geometry of the search models, thus they differ from the standard ones in the case of acting vector field. The study also shows that it is not necessarily efficient to orient the search model in each case such that starting leg is oriented only with or against acting river-type perturbation, as it is routinely assumed. As the main criterium is time under complete search condition \eqref{complete_search}, combining the standard models with the time-optimal paths or creating the new models based solely on the time-optimal paths in the presence of the perturbation may give the potential benefit, i.e higher efficiency of searching. The remodeling which we considered becomes more significant when the ratio $\frac{|W|}{|u|}$ increases. Thus, in the context of some real applications let us note that the current technology enables to implement the models based on the time-optimal paths, for instance in route planning and monitoring referring to the dron aerial survey and patrolling fixed zone, robotic sailing, underwater slow-speed gliding, weather routing combined with the numerical weather prediction models. 

Application to the search models can also be transformed into the non-Euclidean geometric structures. The idea of combining the theory of search with time-optimal paths, represented here by the Randers geodesics, can be followed and developed then. The examples of the extended and generalized theoretical problems arising from our study have been proposed in \S \ref{gen}. Regarding implementations, the applied simplifications admitting only constant perturbation can be preliminarily optimized by considering the stationary current with the use of time-optimal paths as presented in the above examples. Note that the complete solution to Zermelo's problem refers to the perturbation given as a function of position and time, so the standard search models may become inefficient with respect to the criterium of time. This fact gives meaningful opportunity to apply more advanced than the standard search models due to the essential time reduction, in the scenarios with the presence of acting perturbations.


\bigskip

\

\noindent
\textsc{Faculty of Mathematics and Computer Science, Jagiellonian University,\\
6, Prof. S. Łojasiewicza, 30 - 348, Kraków, Poland}\\
and \\
\textsc{Faculty of Navigation, Gdynia Maritime University,\\
3,  Al. Jana Pawla II, 81-345, Gdynia, Poland}\\

\noindent
\textit{E-mail address:} \texttt{piotr.kopacz@im.uj.edu.pl}

\bibliographystyle{siam}
\bibliography{pk_bjga}

\end{document}